\documentclass[english,a4paper,12pt]{article}

\usepackage[intlimits,sumlimits,namelimits]{amsmath}
\usepackage{babel,amsfonts,amsthm,epsfig,vmargin}

\setpapersize{A4}
\setmarginsrb{30mm}{14mm}{20mm}{15mm}{12pt}{9mm}{0pt}{5mm}

\newcommand{\R}{\mathbb{R}}

\newcommand{\E}{\mathcal{E}}
\newcommand{\I}{\mathcal{I}}
\newcommand{\inprod}[2]{#1 \cdot #2}
\newcommand{\der}[1]{\frac{d}{d #1}}

\newcommand{\osader}[1]{\frac{\partial}{\partial #1}}

\theoremstyle{definition}
\newtheorem{lem}{Lemma}[section]

\theoremstyle{remark}

\title{An Efficient Geometric Integrator for Thermostatted Anti-/Ferromagnetic Models}
\author{Teijo Arponen \\
  Institute of Mathematics \\
  Helsinki University of Technology \\
  Finland \\  \and
 Ben Leimkuhler \\
  Department of Mathematics \\
  University of Leicester \\
  U.K.}  
\date{\today}

\begin{document}
\maketitle

\begin{abstract}
  (Anti)-/ferromagnetic Heisenberg spin models arise from discretization of Landau-Lifshitz models in micromagnetic modelling.    In many applications it is essential to study the behavior of the system at a fixed temperature.   A formulation for thermostatted spin dynamics was given by Bulgac and Kusnetsov \cite{bu-ku90:CEA}, which incorporates a complicated nonlinear dissipation/driving term while preserving spin length.   It is essential to properly model this term in simulation, and simplified schemes give poor numerical performance, e.g. requiring an excessively small timestep for stable integration.
 In this paper we present an efficient, structure-preserving method for thermostatted spin dynamics.
\end{abstract}
{\bf Keywords:} Heisenberg ferromagnet, micromagnetics, spin dynamics,
Landau-Lifschitz equation, Gilbert damping, thermostats, constant
temperature, domain walls, geometric integrator, reversible method

\section{Introduction}

In recent years geometric integrators have become ubiquitous for numerical treatment of differential equations.  By a geometric integrator is meant a
numerical method that preserves some known structure of the continuous
flow.   Geometric integrators are particularly important for long term simulations, as used in molecular sampling or celestial mechanics.  In this paper we consider the application of geometric integration principles for the types of spin dynamics systems which arise frequently in modelling of ferromagnets and anti-ferromagnets.   Efficient Lie-Poisson schemes for classical spin dynamics described by the Landau-Lifshitz (LL) equation were studied in \cite{fr-hu-le97:GIC}, and related multisymplectic schemes in \cite{fr03:GST}.     Here we develop and test a geometric integrator for a semi-discrete Landau-Lifshitz-Gilbert
(LLG) equation which includes a nonlinear dissipative term.   This dissipative system forms the
foundation for a more complicated thermostatted model, following the approach of Bulgac and Kusnetsov \cite{bu-ku90:CEA,an-etal96:SDM}.    We design an effective splitting technique for the full coupled system.

LL and LLG are currently a very active topic of research.  Other
approaches to them can be found in
\cite{fr03:GST,wa-ga-e01:GSP,e-wa00:NML,fr-hu-le97:GIC,pr01:CM,sl-ci03:NSN,la-etal99:SDS,le-ni03:GIS},
who also provide further references.
However none of these consider a thermostatted version.

Simulation with the thermostatted version shows fascinating global
behavior: the system first arranges into patterns (spin domains)
with slowly moving domain walls,
then goes into a quasi-chaotic state and quickly rearranges itself into
completely new spin domains.  This kind of transition would not be
possible with local interactions only.  Here the thermostatting
variable is defined in such a way that it has a global character.

A number of recent articles have focussed on the geometric integration of molecular systems in the canonical ensemble \cite{bo-la-le99:NPM,le02:SND,ba-la-le03:GDD,la-le03:GDT}.  In these articles, the aim has been to start from a Hamiltonian formulation for thermostatted molecular simulation and then to provide a suitable symplectic integrator.   The starting point is usually Nos{\'e} dynamics, although generalizations are possible.  

Since the constant energy Heisenberg spin system is Lie-Poisson,  it is natural to seek a Lie-Poisson system to model the action of the thermostat.   While it is possible (with some additional complication, due to the presence of constraints) to develop such a model for the thermostatted Heisenberg model, based on the ideas in \cite{la-le03:GDT}, it is much different in character from the corresponding molecular dynamics models (see the appendix).   In particular, this approach appears to require introduction of many thermostatting variables which act differently on each spin vector of the system.  In the context of magnetic models, this approach therefore sacrifices an important feature of Nos{\'e} molecular dynamics: the apparent compatibility between the thermostatted quasi-dynamics and the microcanonical dynamics.  (Even though Nos{\'e} dynamics is typically only validated based on a phase-sampling correspondence, there is widespread agreement that the thermostatted dynamics is relevant for modelling dynamics of an appropriate extended system which is not too different in character from the microcanonical version.)       Moreover, the Lie-Poisson thermostats add additional complexity in the form of a relatively complex bath model.    

Given these complications, we believe the best available starting point for geometric integration of thermostatted spin dynamics is the alternative framework of Bulgac and Kusnetsov, based loosely on Nos{\'e}-Hoover (NH) dynamics.   Like NH molecular dynamics, these formulations sacrifice Hamiltonian structure, while retaining a reversing symmetry.  It is unclear the extent to which this loss of structure affects the stability of methods and the ultimate resolution of macroscopic features of the spin model.  Although in molecular dynamics it is known that the reversible-only methods are often inferior to their symplectic counterparts \cite{la-st00:CTP}, it is also well established that NH-type methods are far superior to methods that are neither symplectic nor reversible.

The rest of the paper is organized as follows: in Section
\ref{sec:split} we review splitting methods and apply them to our
models.  In Sections \ref{sec:LL_orig}-\ref{sec:thermo} we present the
models and methods in detail.  In Section \ref{sec:numer} we present
numerical results.  Finally in Section \ref{sec:discuss} we present
some conclusions and discussion.

\subsection{Background: Review of splitting methods}
\label{sec:split}

The reader is referred to \cite{le-re04:GIHM,mc-qu02:SM} for a detailed
discussion of splitting methods.  To briefly describe their basic
construction, consider a differential equation $\dot{u} =f(u)$, with
flow map $\Phi_{\tau,f}$.  If $f=f_1+f_2$, we have $\Phi_{\tau,f} =
\Phi_{\tau,f_1}\circ\Phi_{\tau,f_2} + O(\tau^2)$.  If the flows on
vector fields $f$, $f_1$, and $f_2$ share a first integral, then the
composed map will preserve it as well.  In this way, geometric
integrators can be developed to preserve general classes of Lie
groups.  If the vector field is time-reversible, i.e.  $f(Ru)=-Rf(u)$
for some linear involution $R$, then the symmetric concatenation or
``Strang Splitting'' $\hat{\Phi}_{\tau,f} =
\Phi_{\frac{1}{2}\tau,f_1}\circ\Phi_{\tau,f_2}\circ
\Phi_{\frac{1}{2}\tau,f_1}$, where $f_1$, $f_2$ are reversible vector
fields, gives a time-reversible map ($R \hat{\Phi}_{\tau, f}^{-1}=
\hat{\Phi}_{\tau,f} \circ R$), which, moreover, provides a
second-order approximation of the solution on a finite time interval.
As an example, if $H=H(q,p)=T(p)+V(q)$, the leapfrog
(St\"ormer/Verlet) integrator results from the concatenation $
\hat{\Phi}_{\tau,H} = \Phi_{\frac{1}{2}\tau,V}\circ\Phi_{\tau,T}\circ
\Phi_{\frac{1}{2}\tau,V}.  $

The construction of splitting methods for various types of flows, and
with various orders of accuracy, is discussed in a number of papers
(see, e.g, \cite{yo90:CHO,sa-ca94:NHP}).  Practical splitting-based
geometric integrators have been constructed by mathematicians,
chemists and physicists for a wide variety of important applications,
including the rigid body, general holonomic constraints, particle
accelerator models, and the solar system.  Vector field splittings
were used in \cite{fr-hu-le97:GIC} to obtain efficient time-reversible
integrators for (undamped) spin systems; it is this fundamental scheme
that we have extended in this paper to treatment of dissipative and
thermostatted systems.

\section{The original Landau-Lifshitz model as a Poisson system}
\label{sec:LL_orig}

There are several versions of the Landau-Lifshitz equation depending on
which forces and fields are
taken into account.  The version we use  here is that of \cite{fa-ta87:HMTS},
discarding the external  and demagnetizing
field.   (Schemes for more general formulations would build on the work presented here.)  The equation can be written in the form:
\begin{equation}
  \label{origLL}
  \osader{t}S = S \times \nabla^2 S + S \times DS,
\end{equation}
where $x\in I\times I \subset\R^2$, $I$ an interval, $S(x,t)$ is a
unit vector in $\R^3$ representing the classical spin at position $x$
and time $t$, and $D$ is a diagonal matrix representing anisotropy.
Clearly $|S|=$ constant in time:
\begin{equation*}
  \osader{t}|S(x,t)|^2 = 2\inprod{S(x,t)}{\osader{t}S(x,t)}=0\quad\forall x,t.
\end{equation*}
Following the usual practice,  we discretize the spatial variable $x$ using second order
central differences on a regular lattice as in \cite{fr-hu-le97:GIC}
so that in the discretized system the unit length property is
conserved.  We then get a Poisson system on a lattice.  Without loss
of generality we may assume the lattice size to be $1$:
\begin{eqnarray*}
  S(x,\cdot) & \mapsto & z_{ij} \\
  \nabla^2 S(x,\cdot) & \mapsto &
  z_{i,j-1} +z_{i,j+1} +z_{i-1,j} +z_{i+1,j} -4 z_{ij},
\end{eqnarray*}
hence \eqref{origLL} becomes
\begin{equation}
  \label{eq:LL}
  \dot{z}_{ij} = z_{ij} \times 
  (z_{i,j-1} +z_{i,j+1} +z_{i-1,j} +z_{i+1,j} -4 z_{ij}) 
  + z_{ij} \times D z_{ij}.
\end{equation}
Note that the $-4 z_{ij}$ term can be dropped out.  Here we have an
$n\times n$ lattice of spins: the variable $z_{ij}$ is on the unit
sphere of $\R^3$ when $i,j\in\{1,\dots,n\}$.  When either $i$ or $j$
index is zero or $n+1$, those represent boundaries.  Except for the
case of periodic boundary conditions, these boundary terms are
different from the spins: they are an artefact of discretization, and
do not have a counterpart in the continuum case \eqref{origLL}.
Especially, they are not necessarily of unit length.  We do not
represent equations of motion to them, hence they are assumed
constants.

By periodic boundary conditions we mean
\begin{equation}
  \label{eq:BC}
  z_{0j}=z_{nj},\quad
  z_{i0}=z_{in}.
\end{equation}
Next we define the Poisson structure matrix.  Let us denote
\[
z:=\left[z_{11}^T \; | \;z_{12}^T \; | \; \,\dots | \; \,z_{1n}^T \; |
  \;\,z_{21}^T \; | \; \,z_{22}^T \; | \; \,\dots \; | \;
  \,z_{nn}^T\right]^T, 
\]
i.e. $z$ is a column vector.  For an arbitrary
$v=:[a,\,b,\,c]^T\in\R^3$ we denote
\[
  \hat{v}:=
  \begin{pmatrix}
    0  & -c & b  \\
    c  & 0  & -a \\
    -b & a  & 0
  \end{pmatrix}, \qquad \hat{v}u = v \times u \quad \forall u.
\]
The Poisson structure matrix is defined as the block diagonal
\begin{equation}
  \label{J}
  J(z) := 
  \begin{pmatrix}
    \hat{z}_{11} \\
    & \hat{z}_{12} \\
    & & \hat{z}_{13} \\
    & & & \ddots \\
    & & & & \hat{z}_{nn}
  \end{pmatrix}.
\end{equation}
Now \eqref{eq:LL} becomes 
\[
  \dot{z} = J(z) \nabla H(z), 
\]
when we choose the Hamiltonian $H$
\begin{equation}
\label{eq:1}
  H := -\frac{1}{2} \left(\sum_{i,j} 
    \sum_{(a,b)\in NN(ij)} \inprod{z_{ij}}{z_{ab}}
    + \sum_{i,j} z_{ij}^T \,D z_{ij} +H_0 \right),
\end{equation}
where NN refers to ``nearest neighbours'':
\[
  NN(ij) = \{ z_{i,j-1},\, z_{i,j+1},\, z_{i-1,j},\, z_{i+1,j} \},
\]
and $H_0$ represents the boundaries.  For example, if we have zero
boundaries ($z_{i0}=0, z_{0j}=0, z_{i,n+1}=0, z_{n+1,j}=0$), then
\[
  H_0 := 0,
\]
while if we have periodic boundary conditions, then
\begin{equation}
  \label{eq:2}
  H_0 := \sum_j \inprod{z_{0j}}{z_{1j}}
  +\sum_i \inprod{z_{i0}}{z_{i1}}.
\end{equation}
We can easily extend this to a model covering both ferromagnet and
antiferromagnet case.
\begin{equation}
  \label{eq:6}
  H :=-j_K \frac12 \left(\sum_{i,j} \sum_{NN} \inprod{z_{ij}}{z_{ab}}
    + \sum_{i,j} z_{ij}^T \,D z_{ij} +H_0 \right),
\end{equation}
where $j_K$ is the so called exchange integral \cite{an68:MMM}, assumed
constant here, as in \cite{an-etal96:SDM}, and
\begin{equation*}
  j_K
  \begin{cases}
    >0 & \text{ for ferro} \\
    <0 & \text{ for antiferro.}
  \end{cases}
\end{equation*}
Hence we have the Poisson system:
\begin{equation}
  \label{eq:poisson}
  \dot{z} = J(z)\nabla H(z),\qquad H \text{ as in }\eqref{eq:6}.
\end{equation}
For an individual spin at the lattice point $(i,j)$ this becomes
\begin{equation}
\label{eq:5}
\begin{split}
  \dot{z}_{ij} &= -j_K z_{ij} \times
  \left( \sum_{(a,b)\in NN(ij)} z_{ab} \right)
  -j_K z_{ij} \times D z_{ij} \\
  &= z_{ij} \times \nabla H(z), 
\end{split}
\end{equation}
in both periodic and non-periodic cases.    
From now on we employ the notation
\[
  \sum_{NN(ij)} z \quad := \sum_{(a,b)\in NN(ij)} z_{ab}.
\]
\begin{lem}
  \label{lem:1}
  Any system of the form $\dot{z}=J(z)v(z)$ with \eqref{J} and $v$
  an arbitrary vector function, conserves the spin lengths in time:
  \begin{equation}
    \label{eq:12}
    |z_{ij}(t)|=|z_{ij}(0)| \quad \forall i,j,t.
  \end{equation}
\end{lem}
\begin{proof}
  \[
  \der{t}|z_{ij}|^2 = 2 \inprod{z_{ij}}{\dot{z}_{ij}}
  = 2 \inprod{z_{ij}}{z_{ij} \times v(z)} \equiv 0.
  \] 
\end{proof}
This gives us useful freedom in modelling.
Next, the anisotropy term $DS$ is approximated by an average:
\begin{equation}
  \label{roberts}
  D z_{ij} \mapsto
  D \frac14 (z_{i,j-1} +z_{i,j+1} +z_{i-1,j} +z_{i+1,j}).
\end{equation}
this is sometimes referred to \cite{fr-hu-le97:GIC} as the  {\em Roberts discretization}.  
Now \eqref{eq:5} becomes
\begin{equation}
  \label{new_aniso}
  \begin{split}
    \dot{z}_{ij} &= -j_K z_{ij} \times 
    M (z_{i,j-1} +z_{i,j+1} +z_{i-1,j} +z_{i+1,j}) \\
    & = z_{ij} \times \nabla H(z), 
  \end{split}
\end{equation}
where $M=I+D/4$ is a diagonal matrix and $H$ is modified 
according to \eqref{roberts}.

\subsubsection*{Numerical method}

As we noted above, \eqref{eq:poisson} is a Lie-Poisson system whose
meaning we recall here: we can define
\begin{equation*}
  \{f,g\}(z) := \inprod{\nabla f(z)}{(J(z)\nabla g(z))}, 
\end{equation*}
which fulfills the Jacobi identity
\begin{equation*}
  \{\{f,g\},h\} + \{\{g,h\},f\} + \{\{h,f\},g\} =0, 
\end{equation*}
hence $\{\cdot,\cdot\}$ is a Poisson bracket and $J$ is a Poisson
structure matrix.  Since $J$ is linear with respect to $z$, this
Poisson structure can be derived from a Lie algebra structure, hence
it is called a Lie-Poisson structure.

For a detailed discussion on how to integrate this, see
\cite{fr-hu-le97:GIC}.  To summarize that paper, the best way to
integrate is to split the vector field in even-odd (or red-black) way:
\begin{equation}
  \dot{z}_{ij} = V_1 + V_2, 
\end{equation}
where
\begin{eqnarray*}
  V_1 &=&
  \begin{cases}
    -j_K z_{ij} \times M \sum_{NN(ij)} z, & i+j \text{ even} \\[2mm]
    0, & i+j \text{ odd},
  \end{cases} \\
  V_2 &=&
  \begin{cases}
    0, & i+j \text{ even} \\[2mm]
    -j_K z_{ij} \times M \sum_{NN(ij)} z, & i+j \text{ odd}.
  \end{cases}
\end{eqnarray*}
Now, both of these flows can be explicitly solved.  For example $V_1$:
for $i+j$ odd $z_{ij}(t)=z_{ij}(0)$.  For $i+j$ even, the sum over
$NN(ij)$ includes only pairs $a,b$ with $a+b$ odd, hence they are
constants (during $V_1$).  Likewise in $V_2$ the sum is a constant.
Denote the integrator of $V_1$ by $\hat{\Phi}_{1,t}$ and
that of $V_2$ by $\hat{\Phi}_{2,t}$.  That is,
\begin{equation*}
  {\Phi}_{1,t} = \exp(t V_1),\quad
  {\Phi}_{2,t} = \exp(t V_2).
\end{equation*}

The implemented integrator is a symmetric composition of these exact
flows:
\begin{equation}
  \label{eq:conserv_GI}
  \hat{\Phi}_t:= {\Phi}_{2,\frac{t}{2}} \circ {\Phi}_{1,t}
  \circ {\Phi}_{2,\frac{t}{2}}.
\end{equation}
This integrator 
\begin{itemize}
\item is time reversible
\item conserves spin lengths
\item in isotropic case ($D=I$) preserves energy
\end{itemize}
since ${\Phi}_{1,t}$ and ${\Phi}_{2,t}$ do.  See also Section
\ref{sec:split}.

\section{Dissipated version}
\label{sec:dissip}

It is customary to add a dissipation term to \eqref{origLL}.   In our
case the corresponding dissipated version is derived from
\eqref{new_aniso} and becomes
\begin{equation}
  \label{eq:LLG}
  \dot{z}_{ij} = z_{ij} \times \nabla H(z) 
    + \alpha z_{ij} \times z_{ij} \times \nabla H(z),
\end{equation}
where $\alpha$ is a dissipation constant and the corresponding term is
known as the Gilbert damping term.

Clearly \eqref{eq:LLG} can be written more compactly
\begin{equation} 
  \dot{z} = J \nabla H + \alpha J^2 \nabla H.
\end{equation}
From lemma \ref{lem:1} it follows that $|z|=1$ everywhere, i.e. the
dissipation does not affect spin lengths.   Let us first look at the
Gilbert damping term more closely through the equation
\begin{equation}
  \label{gilbert}
  \dot{z} = \alpha z \times (z \times B),
\end{equation}
where $z\in\R^3$, and $\alpha\in\R$ and $B\in\R^3$ are constants.
Or, more compactly,
\begin{equation*}
  \dot{z} = \alpha J^2 B.
\end{equation*}
This can be explicitly solved.  
Put
\begin{eqnarray*}
  v &:=& \inprod{z}{B}, \\
  w &:=& z \times B,
\end{eqnarray*}
then \eqref{gilbert} is
\begin{eqnarray}
  \dot{v} &=& \alpha(-C_1 + v^2), \\
  \dot{w} &=& \alpha v w, \\
  \dot{z} &=& \alpha z \times w,
\end{eqnarray}
where $C_1>0$ constant,
\[
  C_1 = |z|^2 |B|^2.
\]
We can solve for $v$ (we have assumed $|z(0)|=1$):
\begin{equation}
  v(t) = -|B|\frac{\E^2\,C_2 -1}{\E^2\,C_2 +1},
\end{equation}
where
\begin{eqnarray*}
  \E &:=& \exp(\alpha |B| t),  \\[2mm]
  C_2 &:=& \frac{|B|-v_0}{|B|+v_0}.
\end{eqnarray*}
Note: if $t\to\infty$, then
\begin{eqnarray*}
  \alpha >0 &\Rightarrow v(t)\to -|B| &\Rightarrow z,B \text{ become antiparallel} \\
  \alpha <0 &\Rightarrow v(t)\to |B| &\Rightarrow z,B \text{ become parallel}.
\end{eqnarray*}
Substituting $v$ we can solve for $w$, which is a scalar function times a
constant vector:
\begin{eqnarray*}
  w(t) &=& f(t) w(0), \\[2mm]
  f(t) &=& \frac{\E(C_2+1)}{\E^2 C_2 +1} \longrightarrow 0 
  \text{  as $t\to\infty$, if } \alpha\neq 0.
\end{eqnarray*}
Substituting $w$ we can solve for $z$:
\begin{eqnarray}
  \label{eq:gilb.z1}
  z(t) &=& \exp(g \hat{w}_0) z(0), \\[2mm]
  \label{eq:gilb.z2}
  &=& \cos(g|w_0|) z_0 + \frac{\sin(g|w_0|)}{|w_0|} w_0 \times z_0,
\end{eqnarray}
where
\begin{eqnarray}
  \label{eq:gilb.z3}
  g \equiv g(t) &:=& -\alpha \int_0^t f(\tau)d\tau 
  = \frac{C^2+1}{|B|C} \left(\arctan C -\arctan(C\E)\right), \\[2mm]
  \label{eq:gilb.z4}
  C &:=& \sqrt{C_2} = \sqrt{\frac{|B|-v_0}{|B|+v_0}}.
\end{eqnarray}
Note that the exp above is a matrix exponential, while the sin and cos
are the usual scalar functions.  Here $\exp(g \hat{w}_0)$ is expanded
as a Magnus series \cite{ha-lu-wa02:GNI}: the direction of $w(t)$ is
constant, hence $g \hat{w}_0$ commutes with its integrals and Magnus
series truncates after the first term.  The evaluation of that term is
by Rodriguez' formula, hence \eqref{eq:gilb.z2}.

Evaluating $g$ numerically was a problem because eventually $v$
approaches $\pm |B|$ (physically this means $z$ becomes
(anti-)parallel to $B$) so the $C$ in $g$ becomes zero.  $g$ itself is
not singular, however this presentation is difficult to evaluate.
We used the following Taylor expansions in the 
implementation: if $||B|-v_0|<0.0001$, 
\[ 
  g = -1 + \E +C^2\left( -\frac23 +\E -\frac13 \E^3 \right) + \mathcal{O}(C^4),
\]
and if $||B|+v_0|<0.0001$, 
\[ 
  g = 1 - \E^{-1} +C^{-2}\left( \frac23 -\E^{-1} +\frac13 \E^{-3} \right) 
  + \mathcal{O}(C^{-4}).
\]

\subsubsection*{A Lyapunov function}
Note that
\begin{equation*}
  |w|^2 = (z \times B)\cdot(z \times B) = -(B \times z)\cdot(z \times B)
  = -B\cdot z \times(z \times B) = \frac{-B}{\alpha}\cdot \dot{z} 
  = -\frac{\dot{v}}{\alpha},
\end{equation*}
hence 
\begin{equation*}
  \der{t} v = -\alpha |z \times B|^2 \le 0,\text{  if }\alpha \ge 0.
\end{equation*}
So $v$ is a Lyapunov function, when $\alpha$ is positive.  If
$H(z):=Cv=C\inprod{z}{B}$, $C$ constant scalar, then
\begin{equation}\label{H.lyap}
  \der{t}H = -\alpha C |z \times B|^2,
\end{equation}
that is, $H$ is a lyapunov function iff sgn$(\alpha C)=1$.  Later
sgn$(C)$ chooses between ferromagnet and antiferromagnet.

\subsubsection*{Several spins}
Now we continue from \eqref{eq:LLG}, which  can be written
\begin{equation}
  \dot{z}_{ij} = z_{ij} \times (-j_K) M \sum_{NN(ij)} z
  + \alpha z_{ij} \times z_{ij} \times (-j_K) M \sum_{NN(ij)} z.
\end{equation}
Recall from the previous discussion that $\alpha<0$ implies $z_{ij}$
tends to become parallel to
\[ 
  (-j_K) M \sum_{NN(ij)} z.
\]
This means, see \eqref{H.lyap}, that if $j_K>0$, then sgn$(\alpha
j_K)=-1$ and energy $H$ is decreasing.  In other words, for a
ferromagnet negative $\alpha$ means energy damping.

To summarize, a ferromagnetic or antiferromagnetic spin system subject to Gilbert damping will uniformly dissipate energy for appropriate choice of the sign of the damping coefficient.   Moreover, a Gilbert-damped system is spin-length conserving.

\subsubsection*{Numerical method}

To integrate, we split the vector field in even-odd way
as in the conservative case (Section \ref{sec:LL_orig})
\begin{equation}
  \dot{z}_{ij} = V_1 + V_2 + V_3 + V_4,
\end{equation}
where
\begin{eqnarray*}
  V_1 &=&
  \begin{cases}
    -j_K z_{ij} \times M \sum_{NN(ij)} z, & i+j \text{ even} \\
    0, & i+j \text{ odd},
  \end{cases} \\
  V_2 &=&
  \begin{cases}
    0, & i+j \text{ even} \\
    -j_K z_{ij} \times M \sum_{NN(ij)} z, & i+j \text{ odd},
  \end{cases} \\
  V_3 &=&
  \begin{cases}
    -j_K \alpha z_{ij} \times z_{ij} \times M \sum_{NN(ij)} z,
       & i+j \text{ even} \\
    0, & i+j \text{ odd},
  \end{cases} \\
  V_4 &=&
  \begin{cases}
    0, & i+j \text{ even} \\
    -j_K \alpha z_{ij} \times z_{ij} \times M \sum_{NN(ij)} z,
       & i+j \text{ odd}.
  \end{cases}
\end{eqnarray*}
Now, all these flows can be explicitly solved.  For example 
$V_1$: for $i+j$ odd $z_{ij}(t)=z_{ij}(0)$.  For $i+j$ even,
the sum over $NN(ij)$ includes only pairs $a,b$ with $a+b$ odd,
hence they are constants (during $V_1$).  Likewise in $V_2$, $V_3$,
and $V_4$ the sums include only constants.

Hence in $V_1$ and $V_2$ we solve
\begin{equation}\label{oneterm1}
  \dot{z}_{ij} = z_{ij} \times B,\qquad B \text{ constant},
\end{equation}
and in $V_3$ and $V_4$ we solve
\begin{equation}\label{oneterm2}
  \dot{z}_{ij} = \alpha z_{ij} \times z_{ij} \times B,
  \qquad B \text{ constant},
\end{equation}
which are solved above.  Note that \eqref{oneterm1}, \eqref{oneterm2}
have different $B$'s.  The implemented integrator is a symmetric
composition of these exact flows:
\begin{equation}
  \label{eq:dissip_GI}
  \hat{\Phi}_t:= {\Phi}_{4,\frac{t}{2}} \circ 
  {\Phi}_{3,\frac{t}{2}} \circ 
  {\Phi}_{2,\frac{t}{2}} \circ
  {\Phi}_{1,t}
  \circ {\Phi}_{2,\frac{t}{2}}
  \circ {\Phi}_{3,\frac{t}{2}}
  \circ {\Phi}_{4,\frac{t}{2}},
\end{equation}
where $\Phi_{i,t}=\exp(t\,V_i)$ are the exact flows.

An important feature of our method is that it dissipates energy when
the flow \eqref{eq:LLG} does.  This can be seen in the following way:
we solve the flows $V_1,\dots,V_4$ exactly, hence every step in the
composition \eqref{eq:dissip_GI} follows the energy of the associated 
vector field exactly.
Now $\hat{\Phi}_{\Delta t}$ is a second order method and it follows
the energy evolution with accuracy $O(\Delta t^3)$.  With small enough
time step $\Delta t$ the error is negligible and our method dissipates
the energy.

\section{Thermostatted version}
\label{sec:thermo}

The motivation behind using thermostats is keeping the system around
some constant average temperature.  This is a reasonable assumption for
example in systems with heat baths.  That is, we allow the energy to
fluctuate.  But, at the same time, we want to keep the spin lengths
constant.  This will be carried out by modifying the dissipation term
introduced in previous sections.

We use the thermostatting term suggested in \cite{an-etal96:SDM}:
choose a parameter $T$ (temperature) and compare the system's energy
to it, allow the damping coefficient $\alpha$ to vary with time:
$\alpha=\alpha(t)$ and
\begin{equation}
  \label{eq:an.etal}
  \dot{\alpha} = - \left(\frac{\kappa}{\mathcal{N}T}\right)^2
  \sum_{ij} \left(\I - kT \nabla_{z_{ij}} \right) \cdot
  \left( z_{ij} \times z_{ij} \times \I \right),
\end{equation}
where $T$ is temperature, and $k$ is Boltzmann's constant which we hereafter take to be 1.  $\mathcal{N}$ is number of degrees of
freedom, that is $\mathcal{N}=3n^2$ since we have an $n\times n$
square lattice, $\kappa$ is coupling strength and typically
$\approx\sqrt{\mathcal{N}}$.  Here we take $\kappa/\mathcal{N}:=1/n$.
$\I$ is equal to:
\begin{equation}
  \label{I}
  \I \equiv \nabla_{z_{ij}} H := -j_K M \sum_{NN(ij)} z.
\end{equation}
The thermostatting variable $\alpha$ has been given the nickname
``global demon'' \cite{an-etal96:SDM}, so called due to its non-local
(hence non-physical) character: it affects all spins simultaneously.

So our thermostatted system is
\begin{eqnarray}
  \label{thermo1}
  \dot{z}_{ij} &=& z_{ij} \times \I
  + \alpha z_{ij} \times z_{ij} \times \I \\
  \label{thermo2}
    \dot{\alpha} &=& - \left(\frac{\kappa}{\mathcal{N}T}\right)^2
      \sum_{ij} \left(\I - T \nabla_{z_{ij}} \right) \cdot
      \left( z_{ij} \times z_{ij} \times \I \right).
\end{eqnarray}

It is possible to show that the ferromagnetic system thermostatted using  \eqref{thermo1},\eqref{thermo2} samples from the canonical ensemble.  This system also conserves spin length.    Finally, one easily demonstrates that these equations are invariant under the simultaneous time-coordinate transformation $t\mapsto -t$, $z\mapsto -z$, $\alpha\mapsto -\alpha$, i.e. the equations are time-reversible.

\subsubsection*{Numerical method}

To integrate, we split the vector field in even-odd way as above, with
the $\dot\alpha$ term:
\begin{equation}
  \begin{pmatrix}
    \dot{z}_{ij} \\[2mm]
    \dot{\alpha}
  \end{pmatrix}
  = V_1 + V_2 + V_3 + V_4 + V_5,
\end{equation}
where $V_1,\,V_2,\,V_3,\,V_4$ as in Section \ref{sec:dissip} (with
$\alpha=\alpha_0=$constant) and
\begin{equation}
  V_5 \leftrightarrow
  \begin{cases}
    z_{ij} &= \text{ constant}, \\[2mm]
    \dot\alpha &= \sum_{ij} \left(
      (\inprod{z_{ij}}{B})^2 -\inprod{B}{B} -2T \inprod{z_{ij}}{B} \right).
  \end{cases}
\end{equation}
Here we have simplified:
\[ 
  B:= \I = -j_K\,M \sum_{NN(ij)} z = \text{ indep. of }z_{ij},
\]
\[
  \nabla_z \cdot (z \times z \times B) = 
  \nabla_z \cdot ((\inprod{z}{B})z - B) = 2\inprod{z}{B},
\]
\[
  B \cdot (z \times z \times B) = (\inprod{z}{B})^2 -\inprod{B}{B}.
\]
 
In $V_5$ all terms are constant so the equation with $V_5$ is trivially solved.
But note that the update step of $\alpha$ is $O(n^2)$.
\begin{equation}
  \label{eq:thermo_GI}
  \hat{\Phi}_t:= {\Phi}_{1,\frac{t}{2}} \circ 
  {\Phi}_{2,\frac{t}{2}} \circ 
  {\Phi}_{3,\frac{t}{2}} \circ 
  {\Phi}_{4,\frac{t}{2}} \circ
  {\Phi}_{5,t}
  \circ {\Phi}_{4,\frac{t}{2}}
  \circ {\Phi}_{3,\frac{t}{2}}
  \circ {\Phi}_{2,\frac{t}{2}}
  \circ {\Phi}_{1,\frac{t}{2}},
\end{equation}
where $\Phi_{i,t}=\exp(t\,V_i)$ are the exact flows.

An important feature of this discretization is that it is
time-reversible with respect to the mapping $z\mapsto -z,\quad
\alpha\mapsto -\alpha,\quad t\mapsto -t$.  This can be seen by
recalling from Section \ref{sec:split} that if $f(Ru)=-Rf(u)$ for some
linear involution $R$, then the Strang splitting gives a
time-reversible map.  Here $u:=(z,\alpha)$ and $Ru:=(-z,-\alpha)$.
Applying the rule four times in a row: first to $\Phi_4$ and $\Phi_5$
in the roles of  $\Phi_{\tau,f_1}$ and $\Phi_{\tau,f_2}$ of Section
\ref{sec:split}, secondly to $\Phi_3$ and $\Phi_4\circ\Phi_5$ in a
similar way, next to $\Phi_2$ and $\Phi_3\circ\Phi_4\circ\Phi_5$ and
finally to $\Phi_1$ and $\Phi_2\circ\Phi_3\circ\Phi_4\circ\Phi_5$, we
get the claim.

\section{Numerical results}
\label{sec:numer}

In all our simulations we used $n=50$, that is, a $50\times 50$
lattice.  We used ferromagnets with anisotropy:
$D=$diag$(1,1,\lambda)$.  This is known as ``easy plane'' or ``easy
axis'' anisotropy, corresponding to $\lambda<1$ or $\lambda>1$,
respectively.

\subsection{Dissipated system}
\label{sec:num.dissip}

\begin{figure}[htb]
  \begin{center}
    \epsfig{file=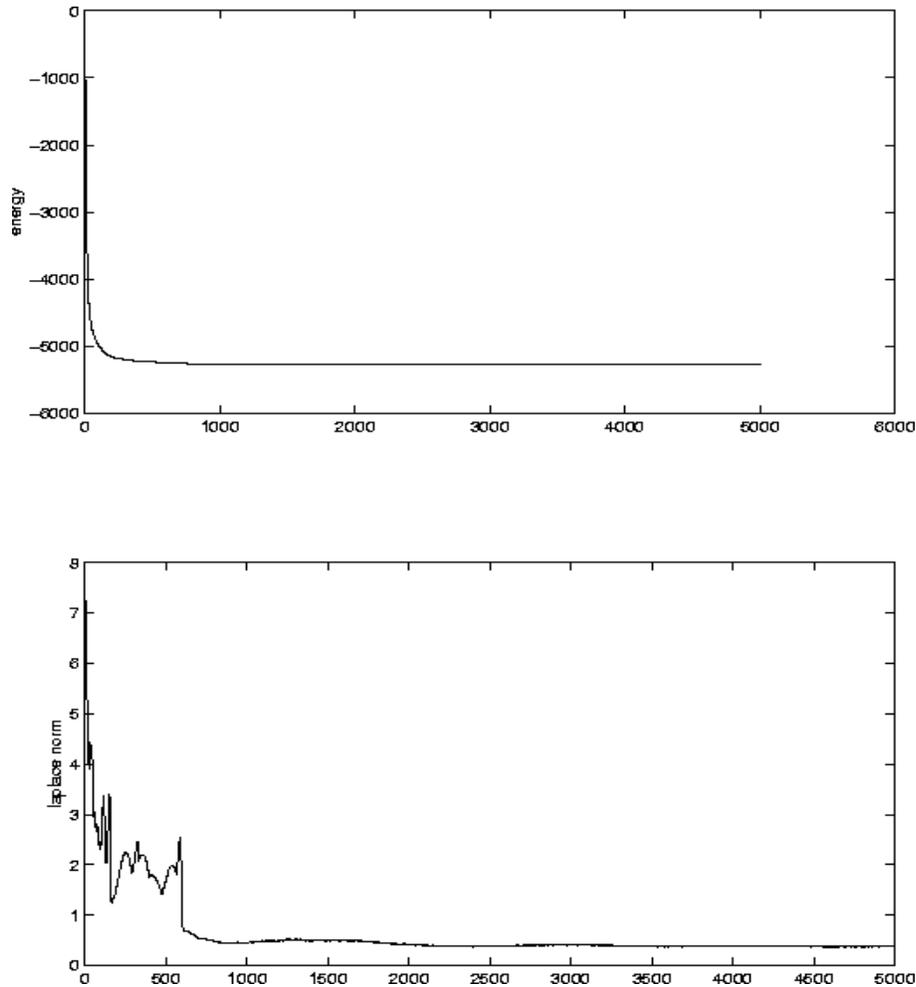,height=13cm,width=12cm}
    \caption{Energy of dissipated system.}
    \label{fig:1}
  \end{center}
\end{figure}

\begin{figure}[htbp]
  \begin{center}
    \begin{tabular}{lr}
      \epsfig{file=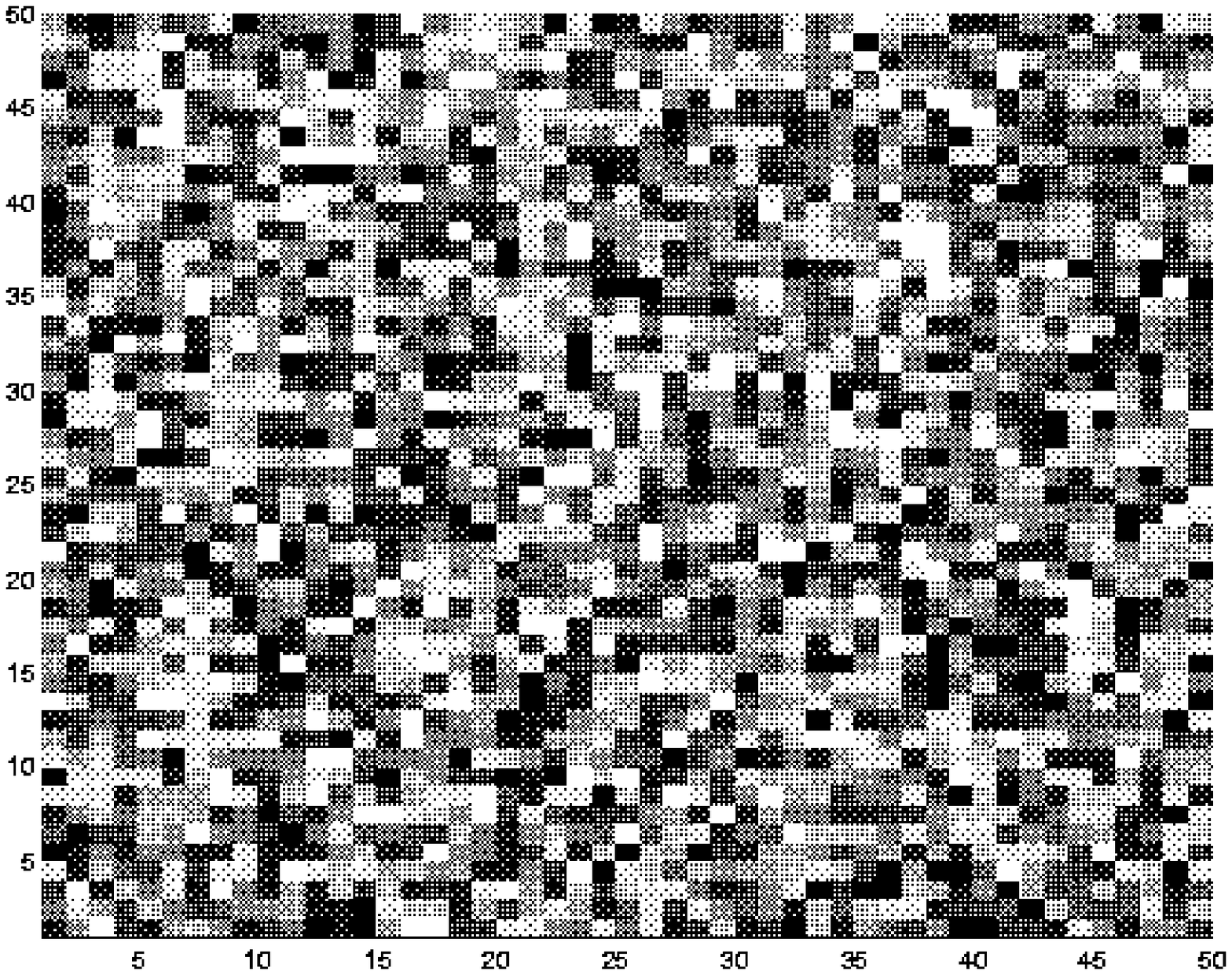,width=7cm}
      & \epsfig{file=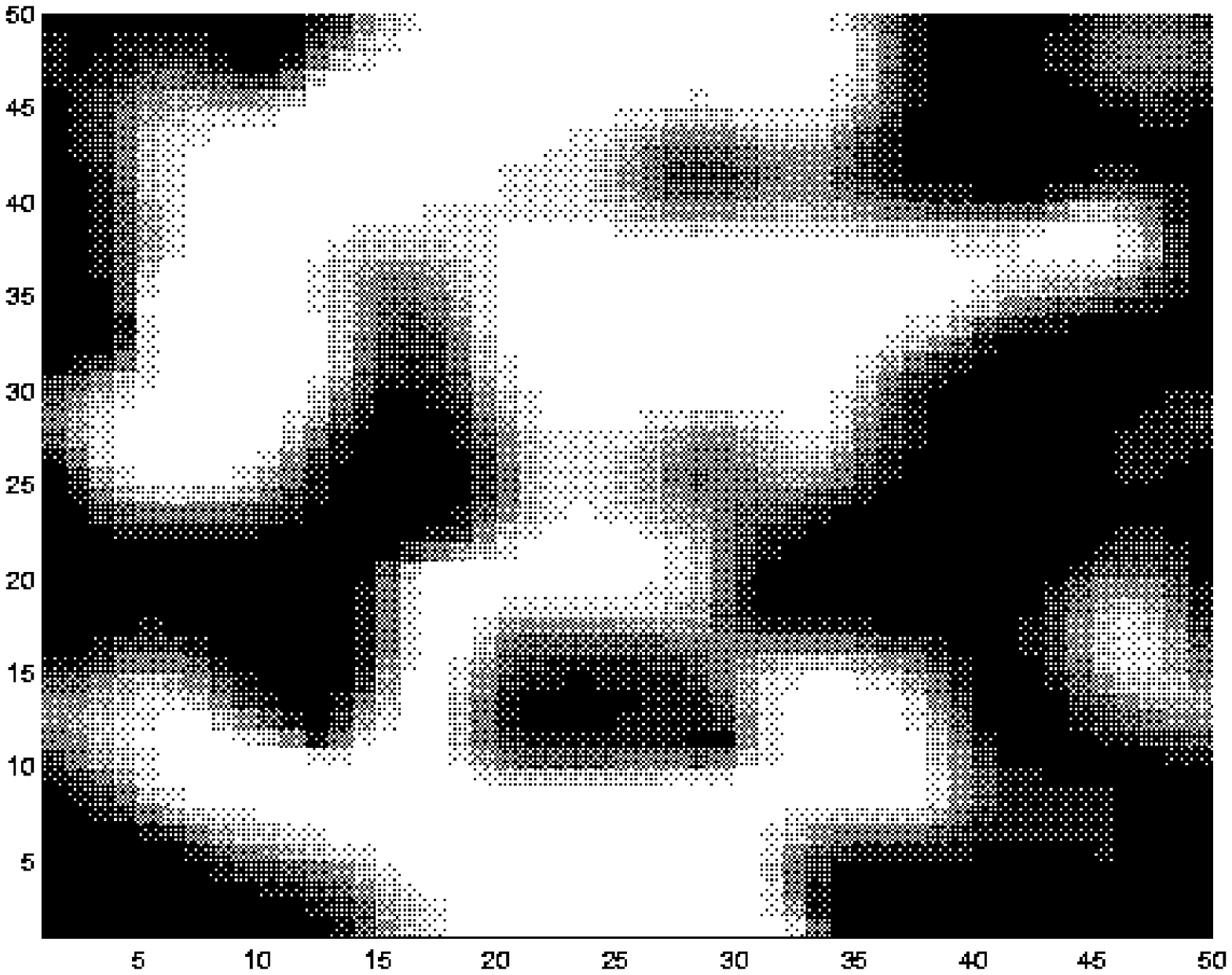,width=7cm} \\
      $ t=0$   &   $ t=0.1$ \\
      \epsfig{file=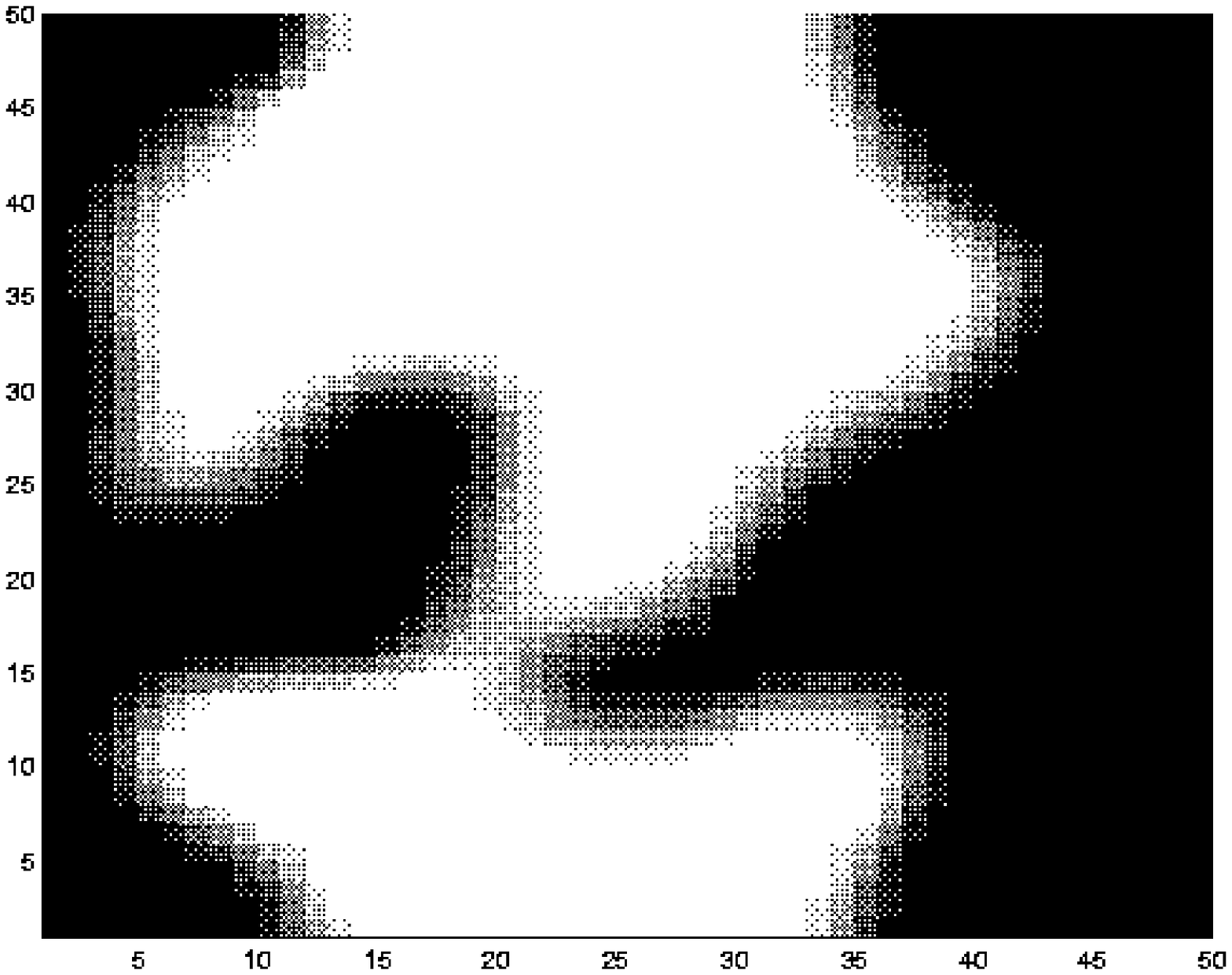,width=7cm}
      & \epsfig{file=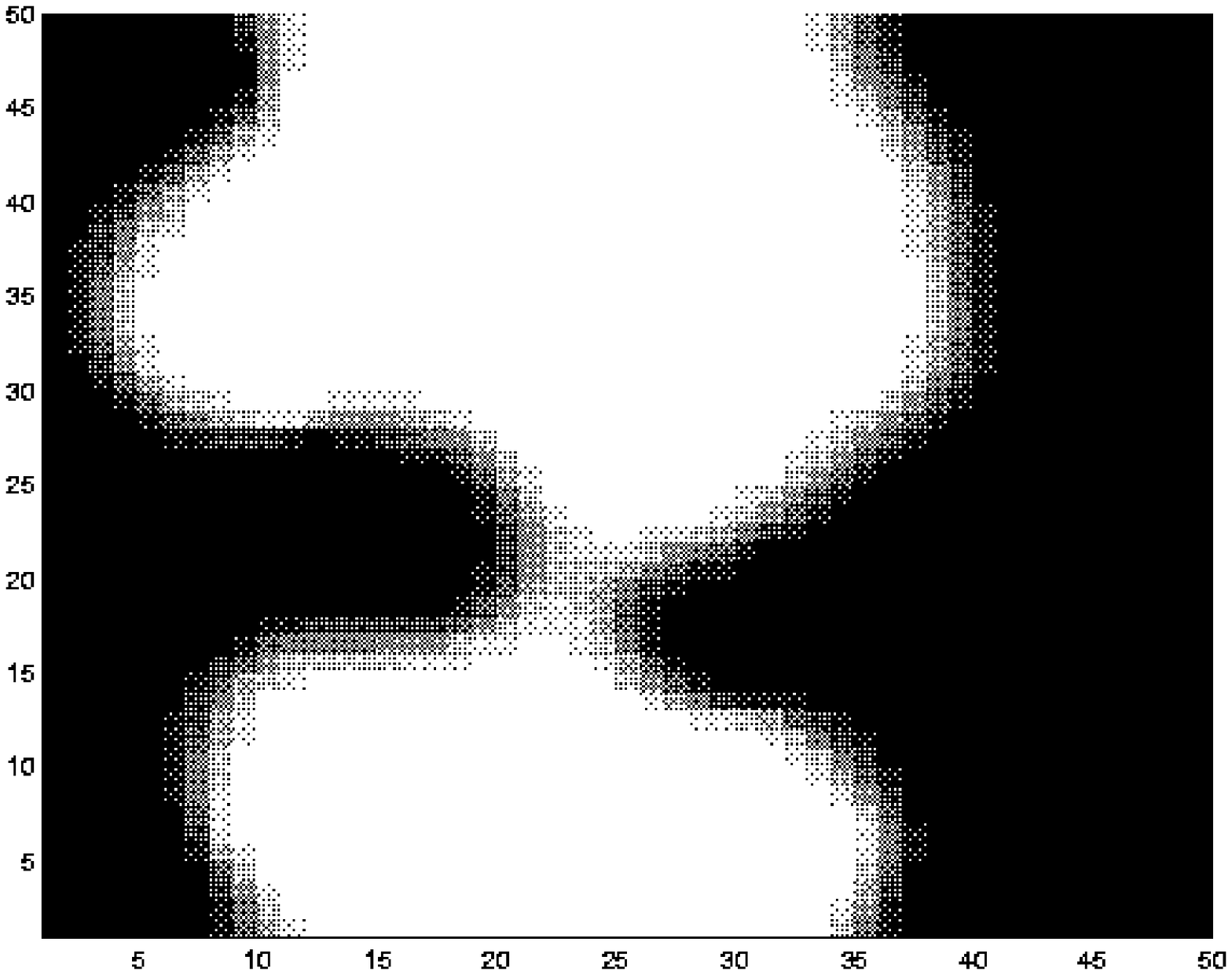,width=7cm} \\
      $ t=0.2$   &   $ t=0.3$ \\
      \epsfig{file=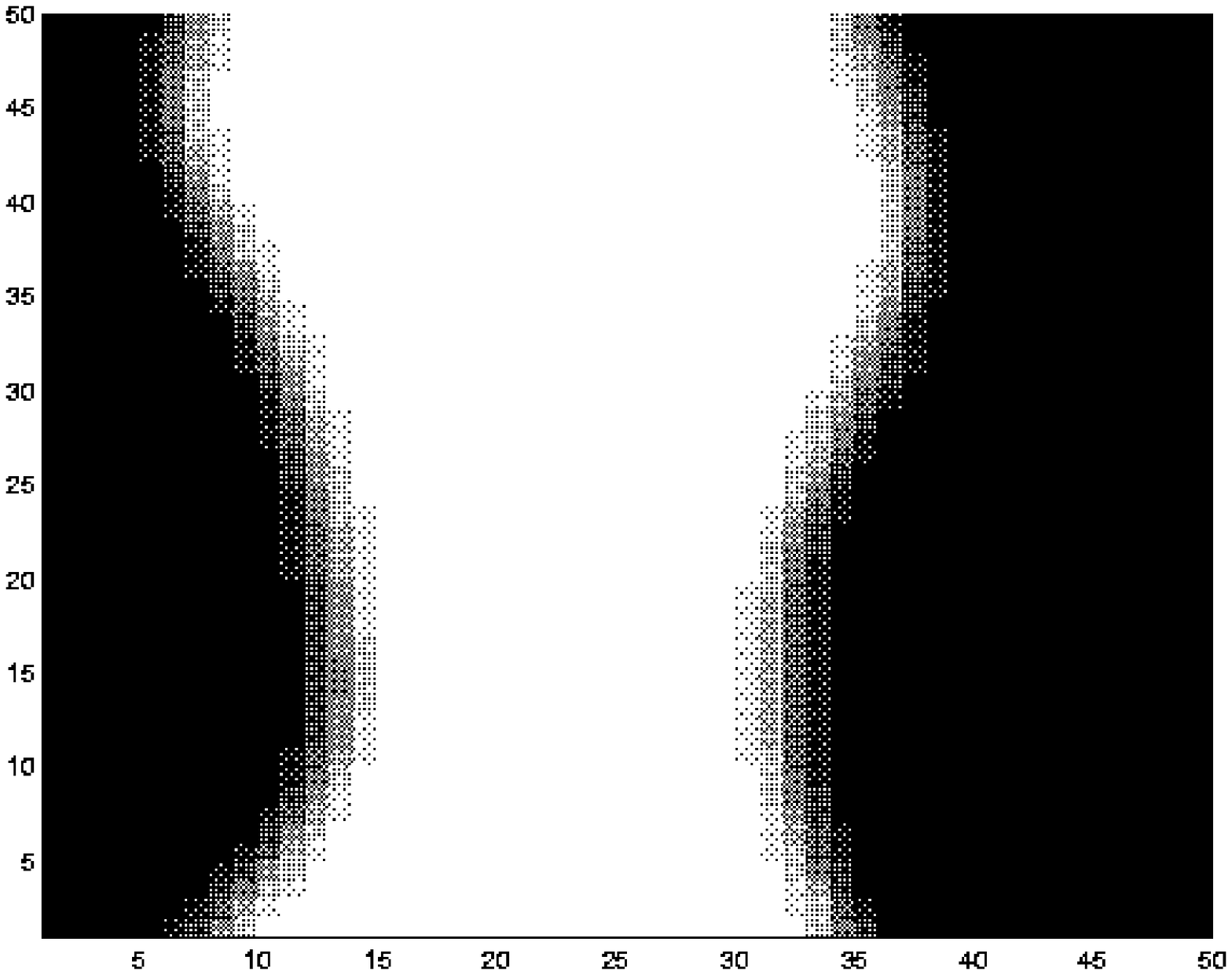,width=7cm}
      & \epsfig{file=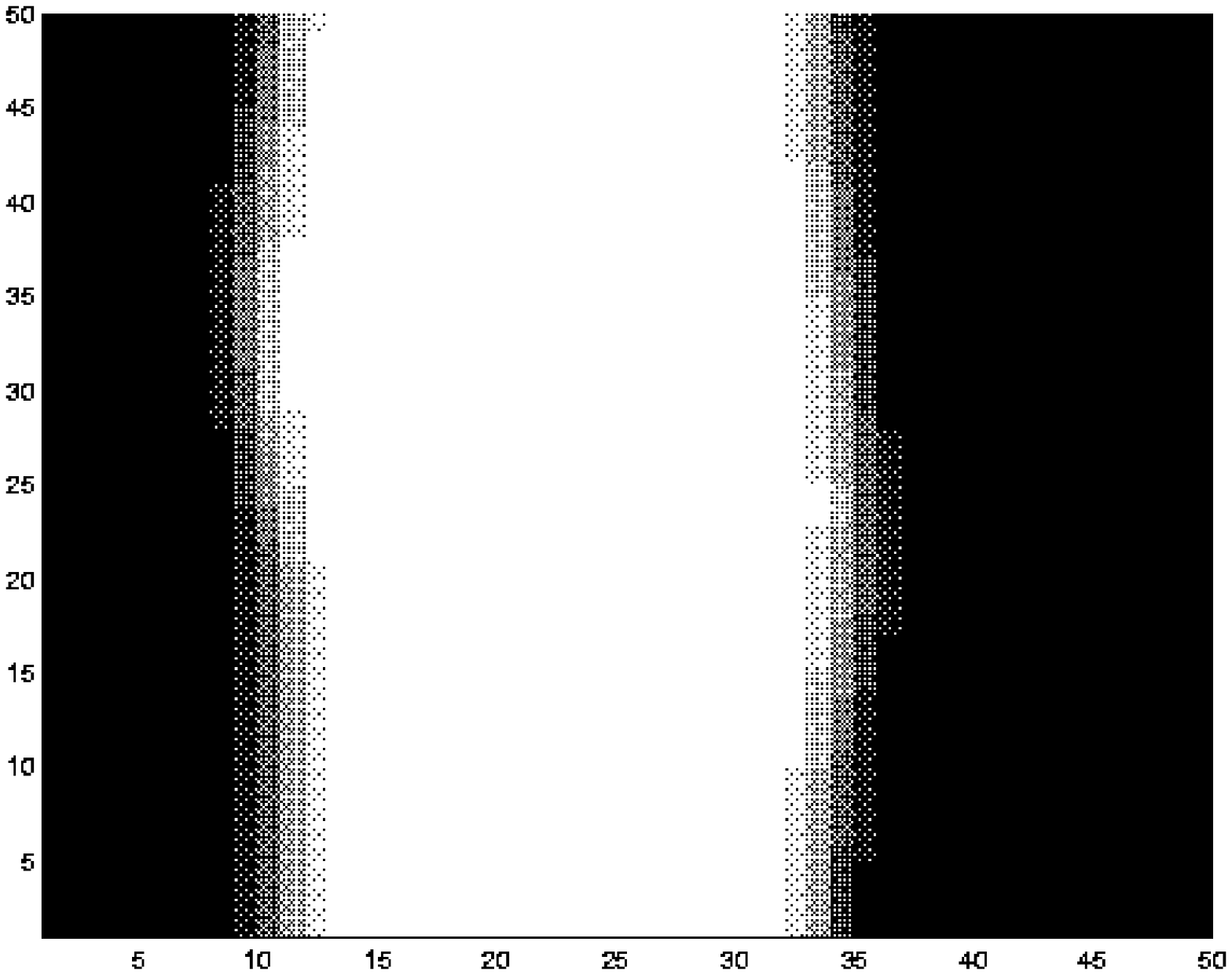,width=7cm} \\
      $ t=1.1$   &   $ t=1.6$
    \end{tabular}
    \caption{Snapshots of $z$-components (black is down spin, white is up spin) of the dissipated system, example 1.}
    \label{fig:2}
  \end{center}
\end{figure}

{\em Example 1.}
In Figure \ref{fig:1} we see evidence of the dissipation of energy.
On the top part is energy in semilog scale, on the bottom part is 
the maximum norm of the discrete Laplacian during each time step.
Here we used periodic boundary conditions, $\lambda=1.1$, $\alpha=-0.5$, and timestep $\Delta t=0.1$.
The initial configuration was random (top left of Figure \ref{fig:2}).

The evolution of the discrete Laplacian is understood as the system
settling down to some formation, and this can be seen in Figure
\ref{fig:2}, which includes snapshots of the same simulation.  The
snapshots describe the $z-$components of the spins.  The order of the
pictures is top row first, from left to right.  The darker a point is,
the lower $z-$component it has.  Black represents spin down, white
spin up.

The result is typical: the dissipated system converges to two bands of
up and down spins.  We also tested zero boundary conditions, then the
dissipated system typically converged to a single band.

\subsection{Thermostatted system}
\label{sec:num.thermo}

\begin{figure}[htb]
  \begin{center}
    \epsfig{file=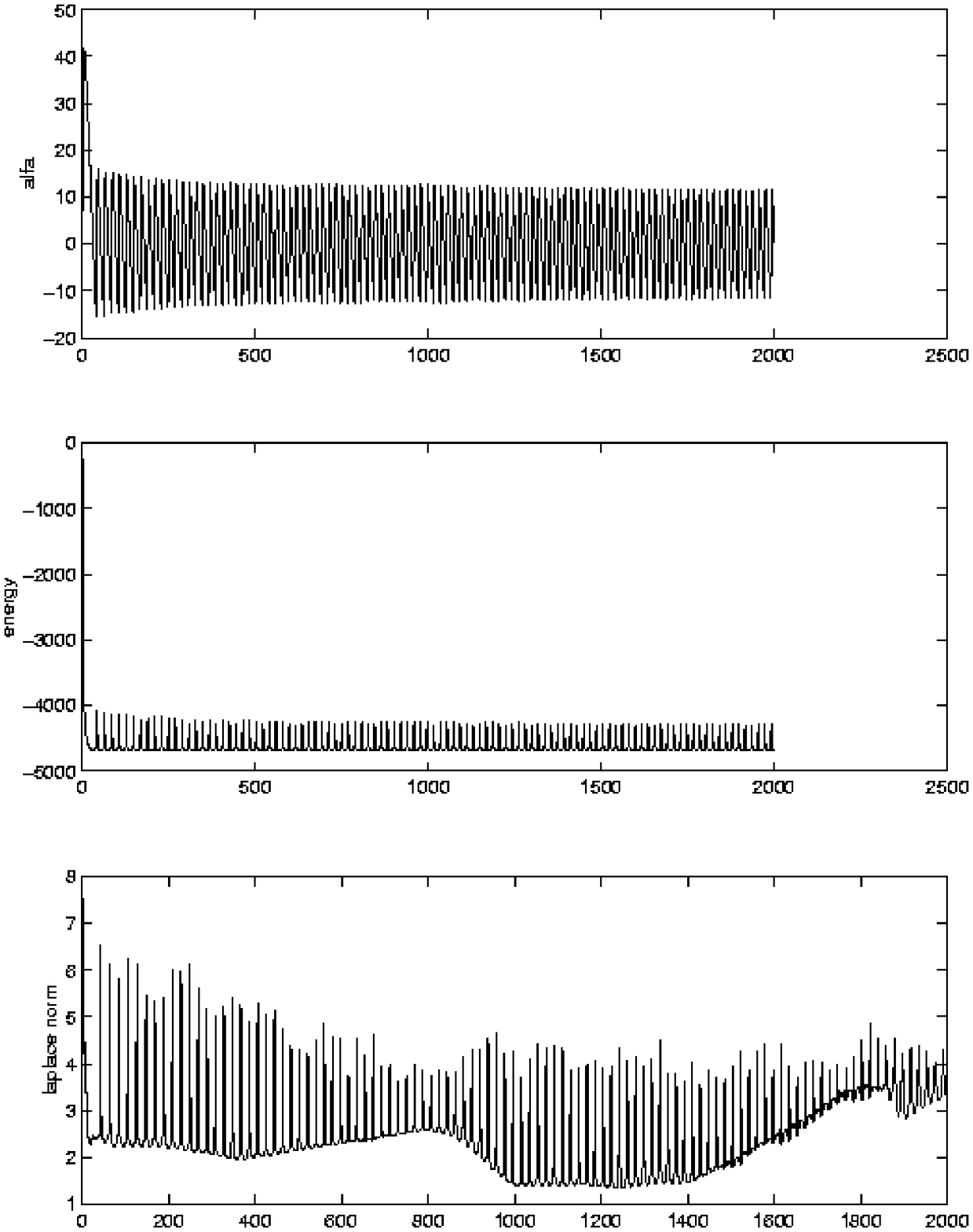,height=13cm,width=12cm}
    \caption{Energy, alpha, and norm of Laplacian of thermostatted system.}
    \label{fig:3}
  \end{center}
\end{figure}

\begin{figure}[htbp]
  \begin{center}
    \begin{tabular}{lr}
      \epsfig{file=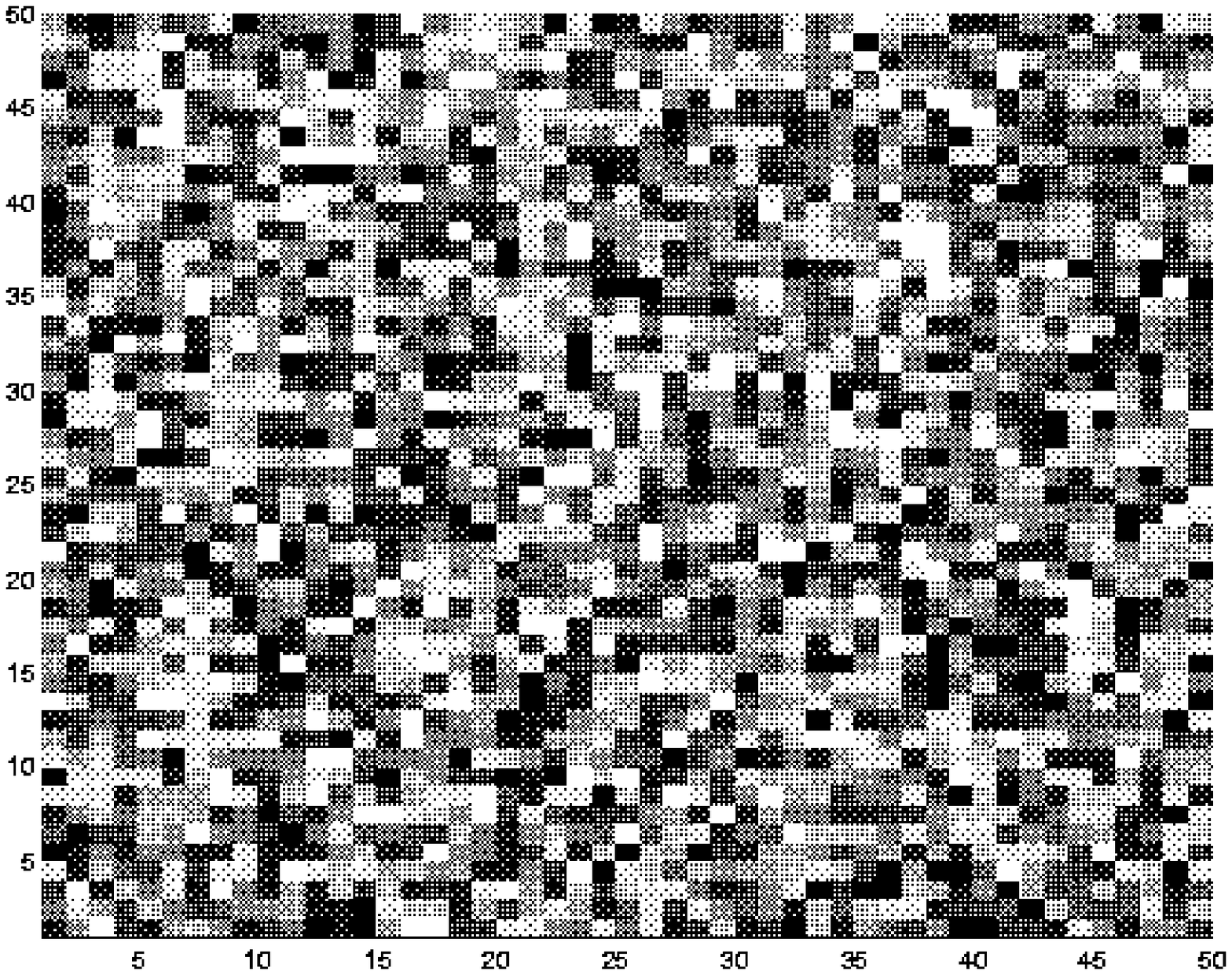,width=7cm}
      & \epsfig{file=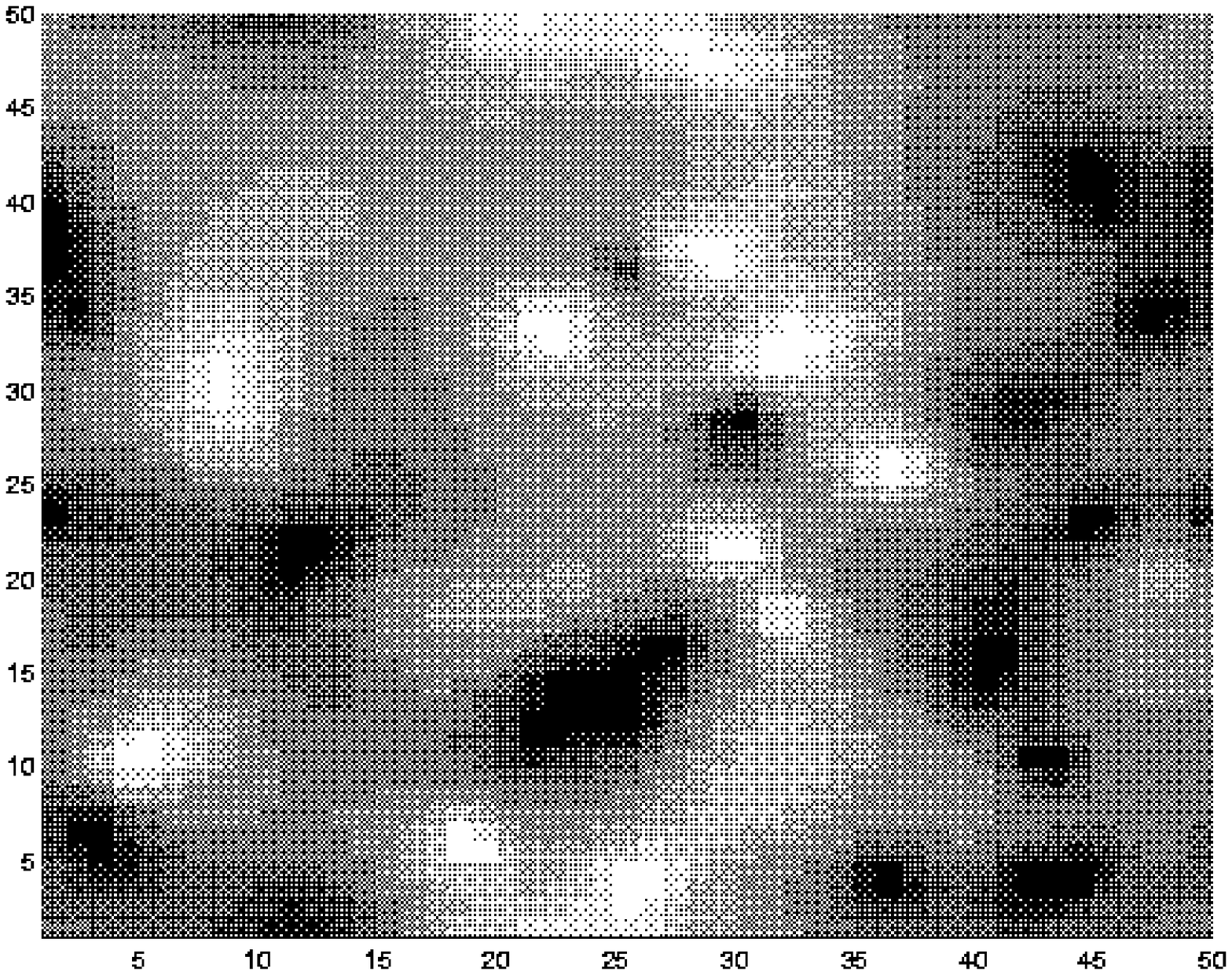,width=7cm} \\
      $ t=0$   &   $ t=0.01$ \\
      \epsfig{file=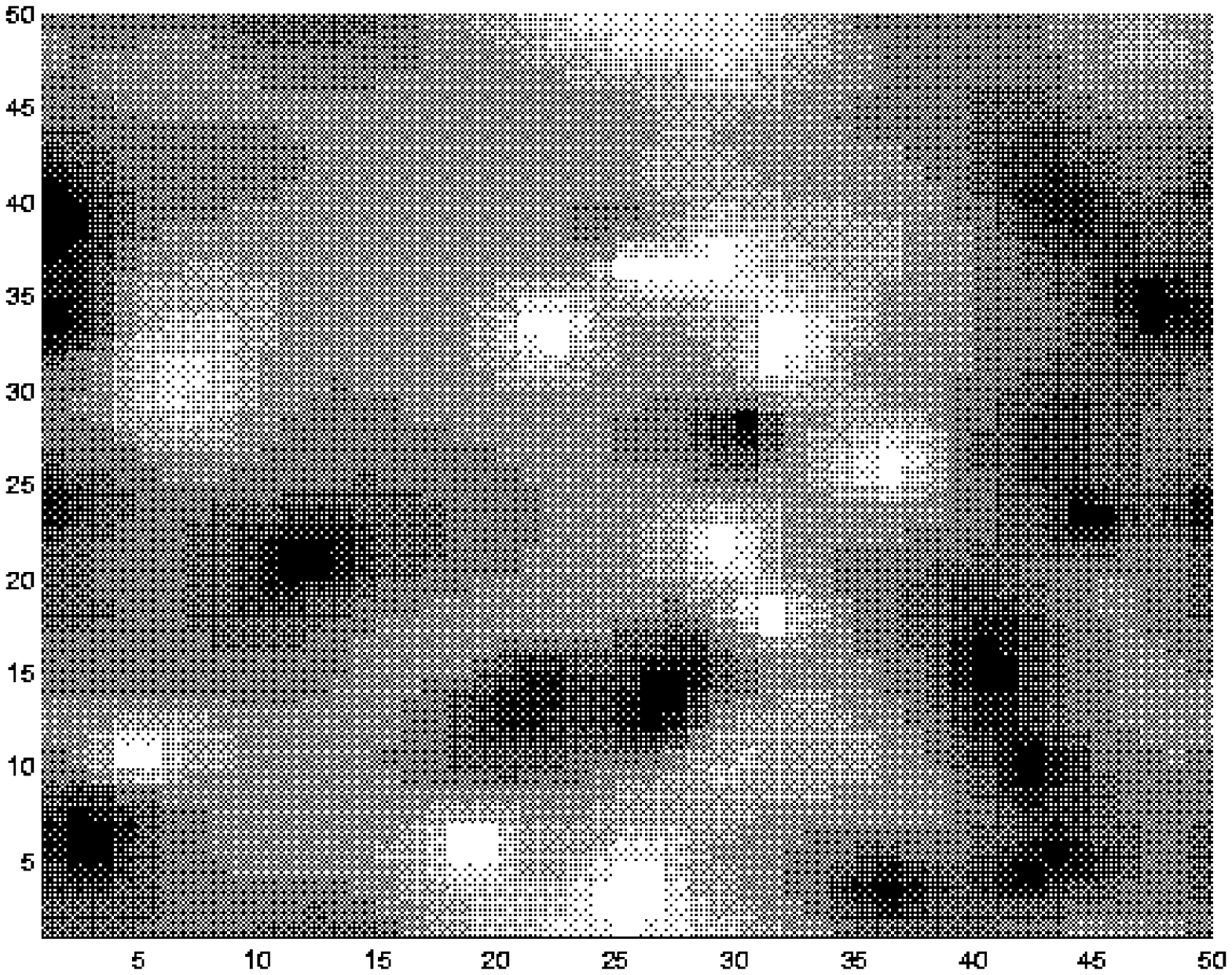,width=7cm}
      & \epsfig{file=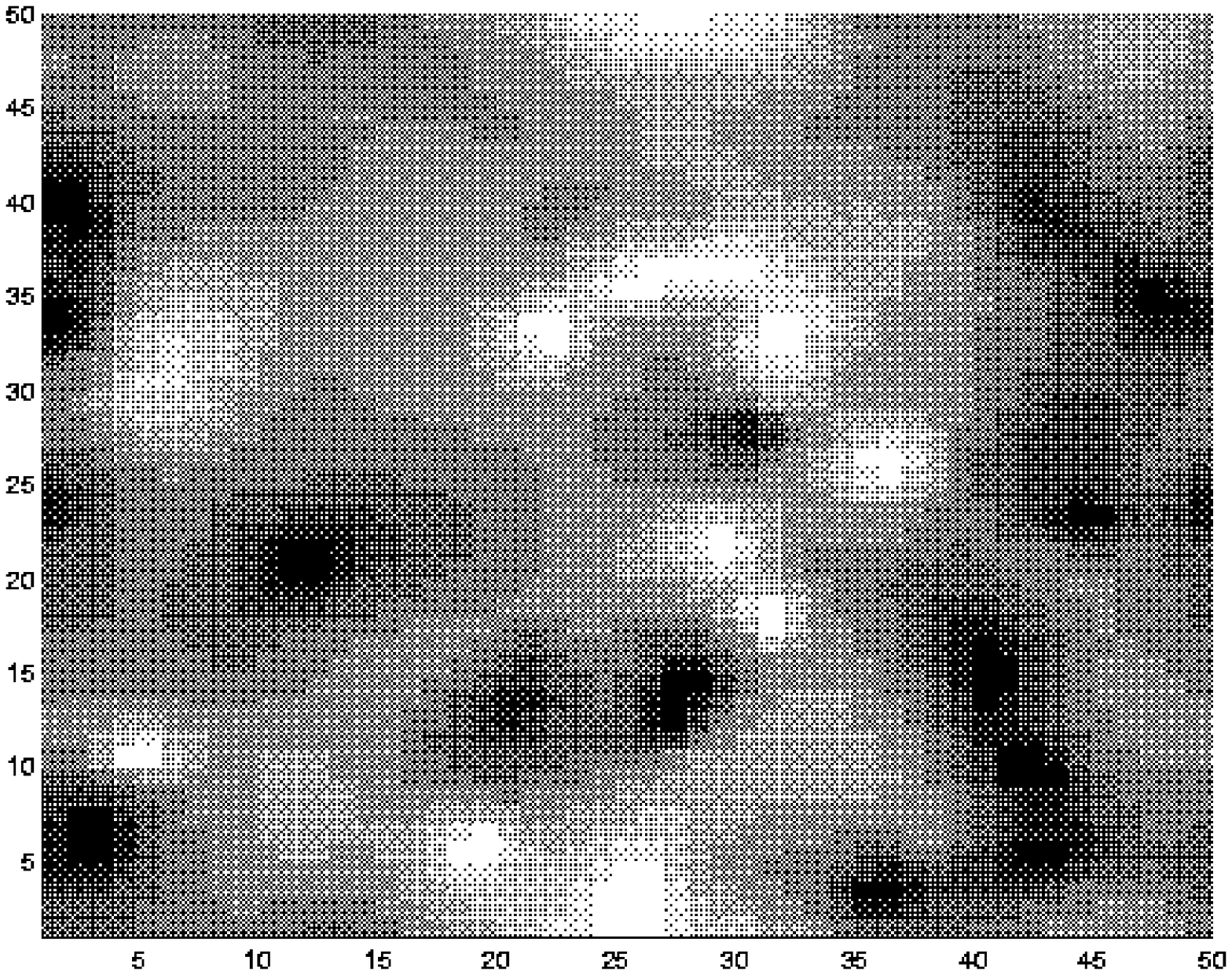,width=7cm} \\
      $ t=0.11$   &   $ t=0.16$ \\
      \epsfig{file=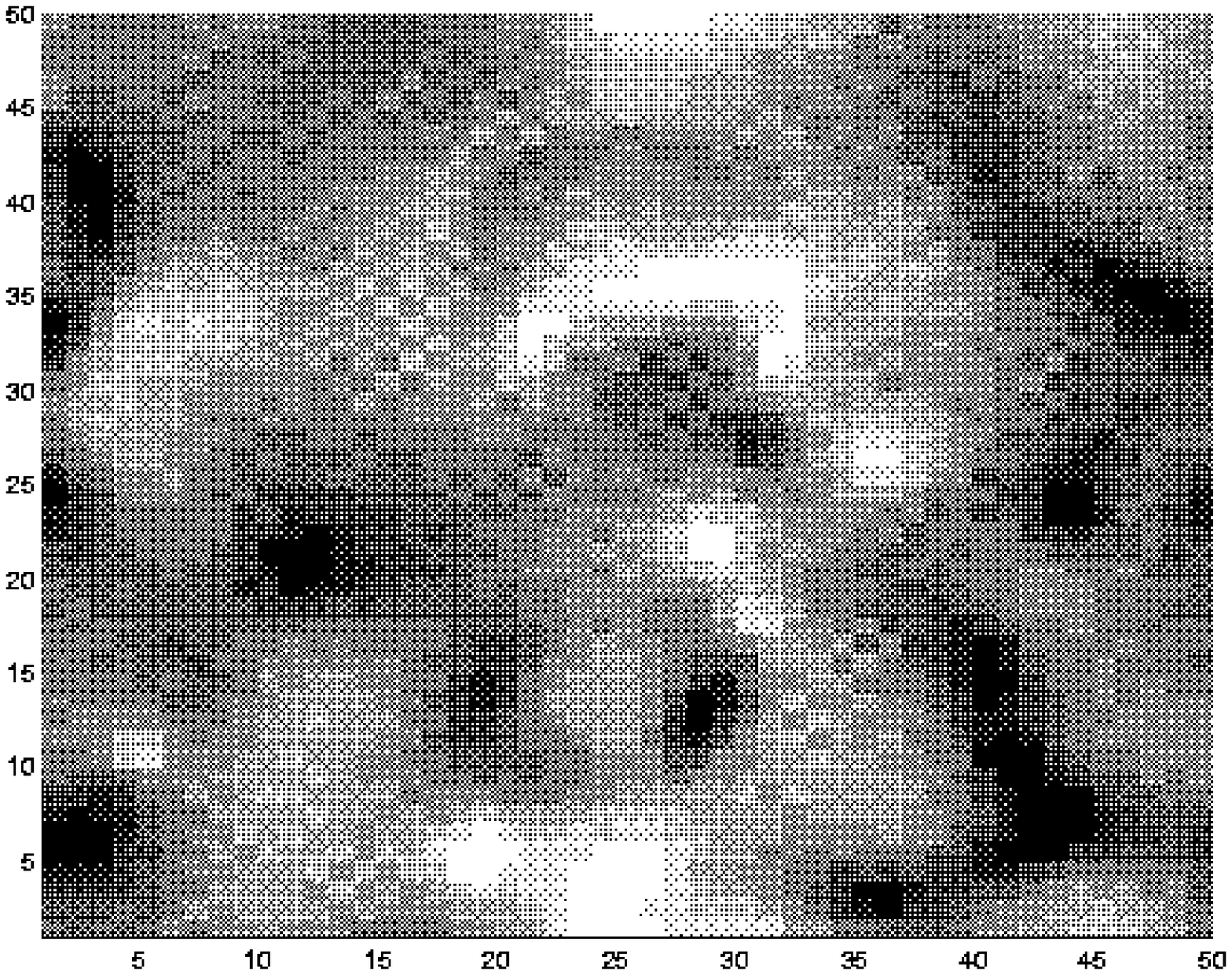,width=7cm}
      & \epsfig{file=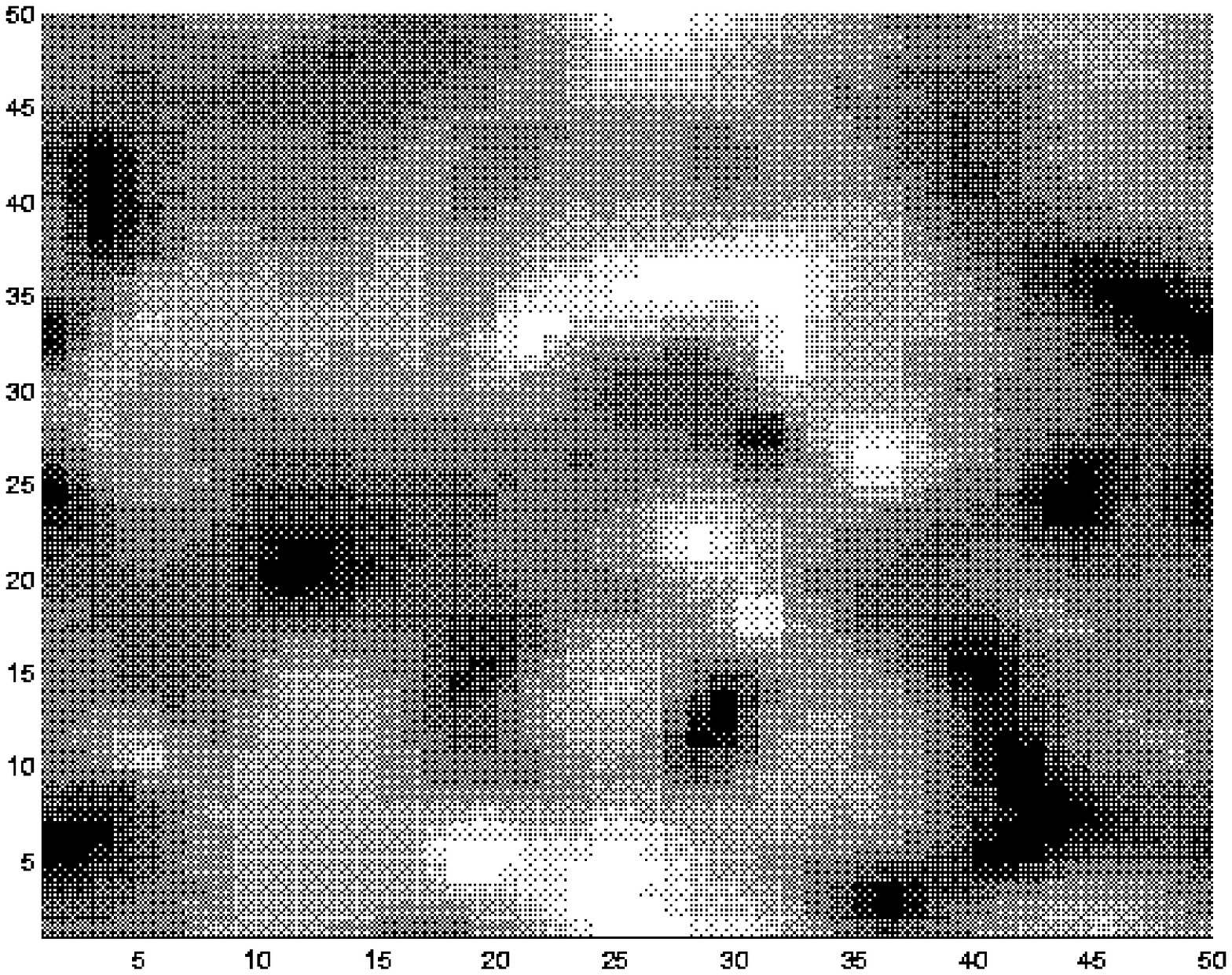,width=7cm} \\
      $ t=0.26$   &   $ t=0.31$
    \end{tabular}
    \caption{Snapshots of $z$-components (black is down spin, white is up spin) of thermostatted system, Example 2.}
    \label{fig:4}
  \end{center}
\end{figure}

{\em Example 2.}  In Figures \ref{fig:3} and \ref{fig:4} is a
thermostatted system with periodic boundaries, $\lambda=0.9$, $T=0.04$
and timestep $\Delta t=0.01$.  The initial condition is random.  In
the top part of Figure \ref{fig:3} is the thermostatting variable
$\alpha$ and in the middle part is the energy, and in the bottom part
the maximum norm of the discrete Laplacian.  After an initial phase
both $\alpha$ and energy settle to an aperiodic oscillatory motion,
$\alpha$ between $-10$ and $+10$, energy between $-4800$ and $-4300$.
We plotted only 2000 steps but the behavior continued similarly for
at least 25000 steps.  In Figure \ref{fig:4} we can see slowly
creeping boundaries; the reader is asked to compare the white areas.
At $t=0.26$ it suddenly looks chaotic, then renders back to the
creeping boundaries.  This can be seen as a kind of stability of the
creeping boundaries.  After 25000 steps there still is slow motion,
the system does not converge to any particular formation.

\begin{figure}[htb]
  \begin{center}
    \epsfig{file=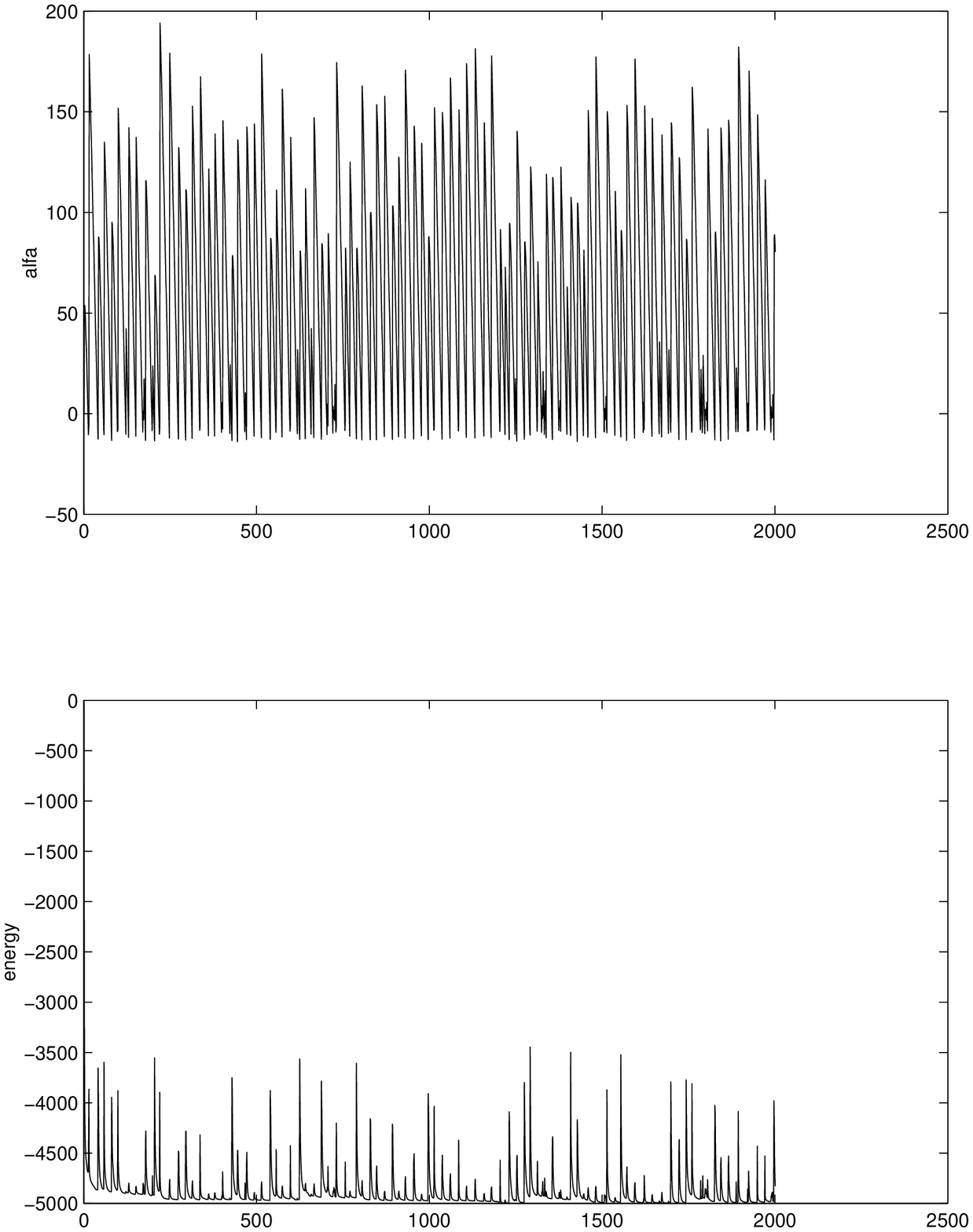,height=13cm,width=12cm}
    \caption{Energy and alpha of ``wandering vortices'', Example 3.}
    \label{fig:5}
  \end{center}
\end{figure}

\begin{figure}[htb]
  \begin{center}
    \epsfig{file=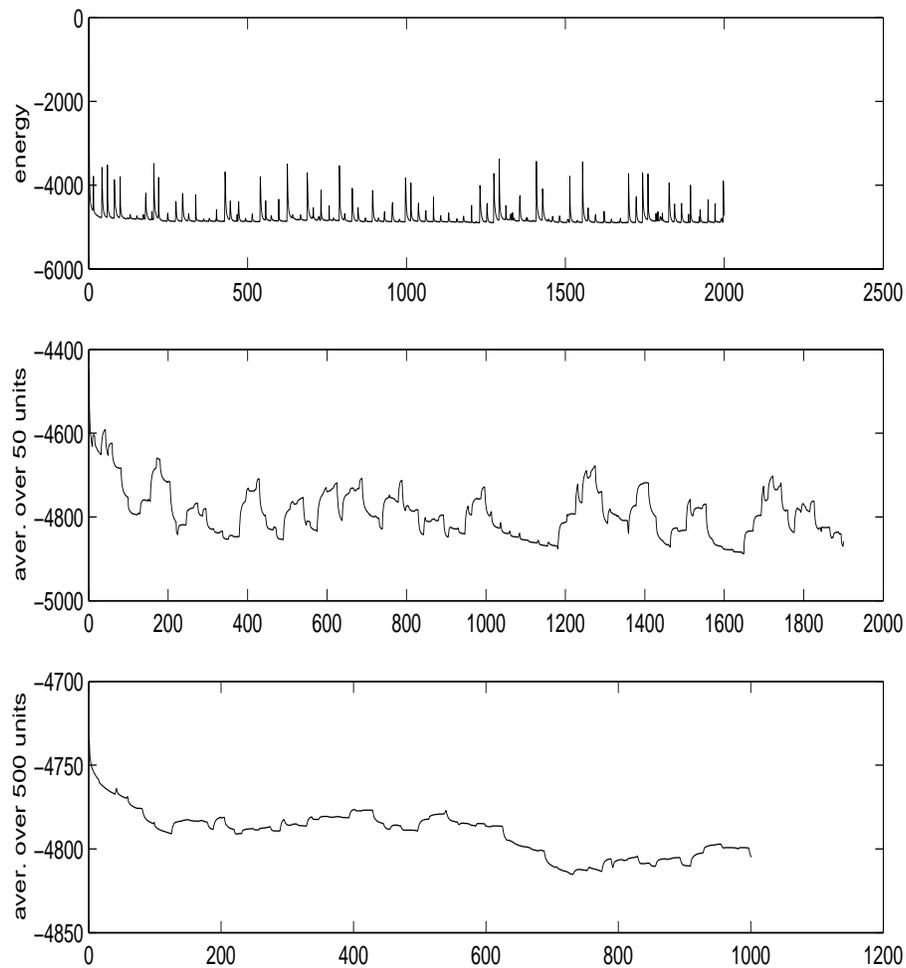,height=13cm,width=12cm}
    \caption{Averages of energy of ``wandering vortices'', Example 3.}
    \label{fig:6}
  \end{center}
\end{figure}

\begin{figure}[htbp]
  \begin{center}
    \begin{tabular}{lr}
      \epsfig{file=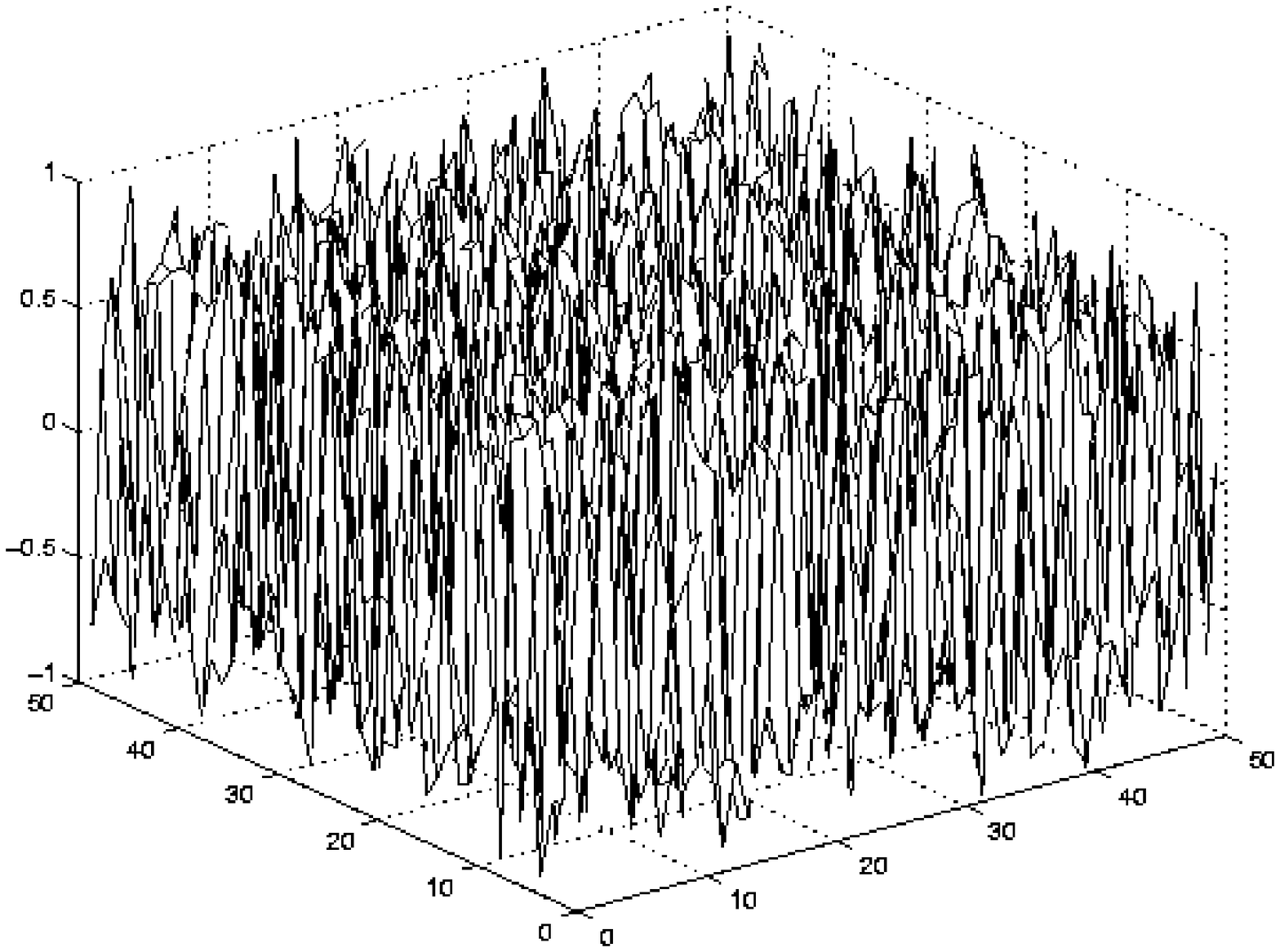,width=7cm}
      & \epsfig{file=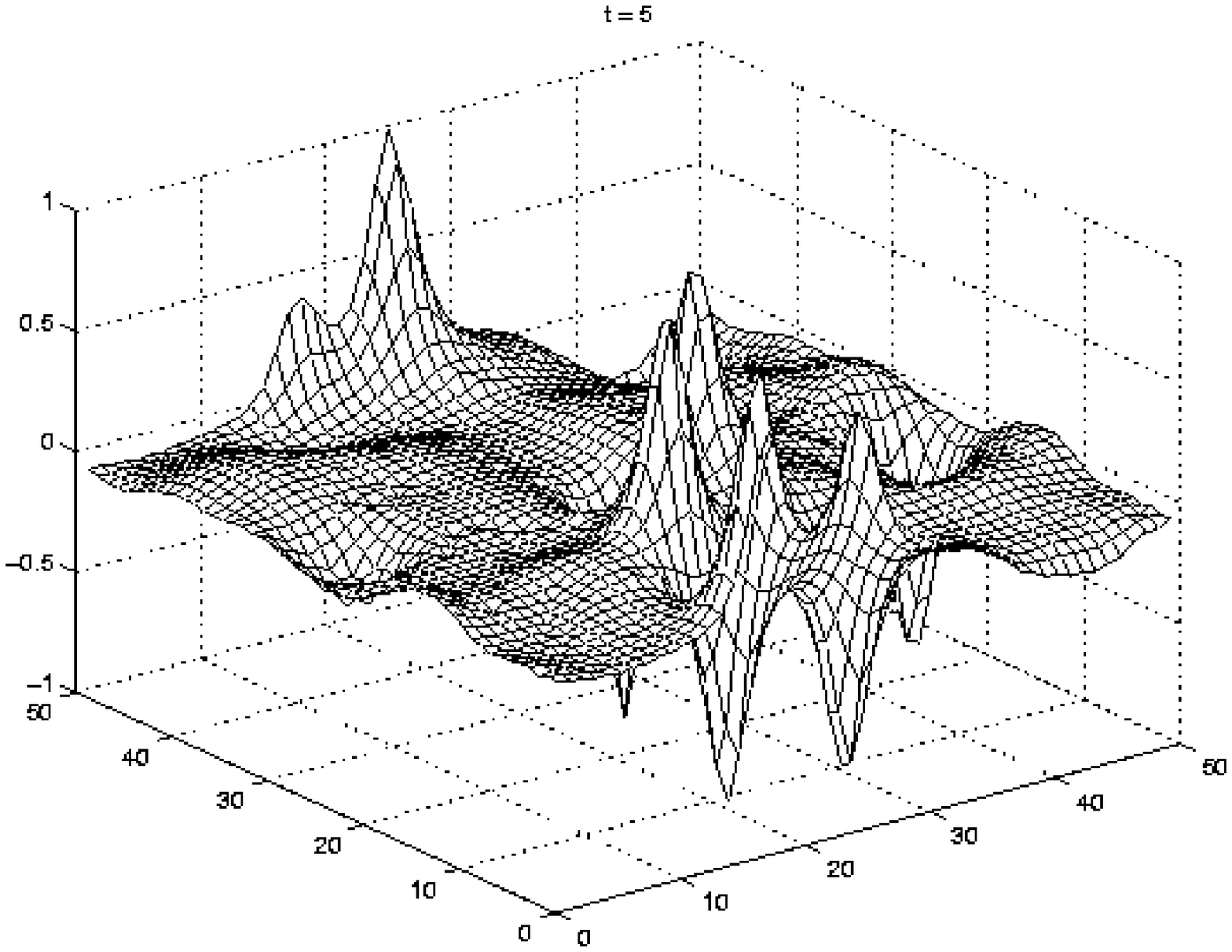,width=7cm} \\
      $ t=0$   &   $ t=5$ \\
      \epsfig{file=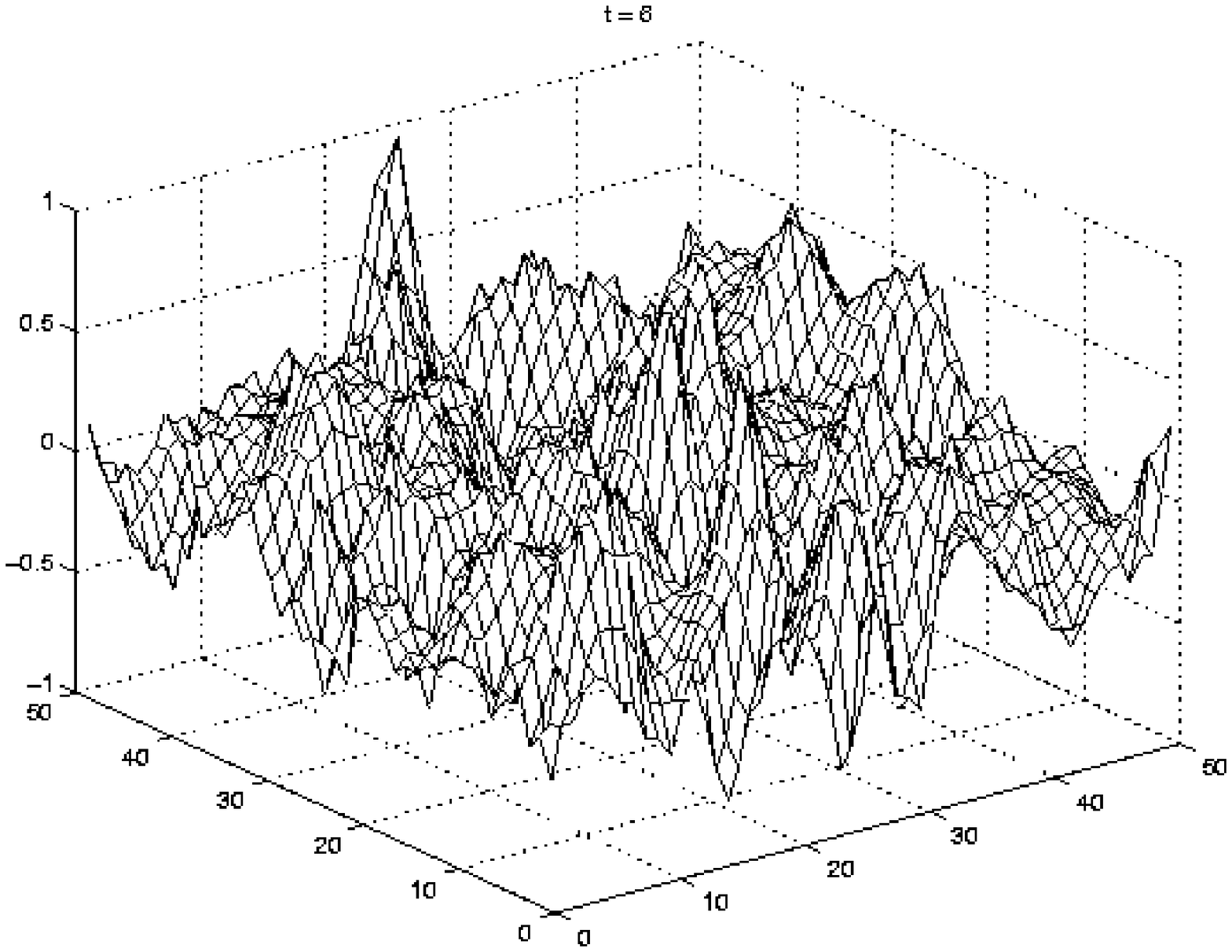,width=7cm}
      & \epsfig{file=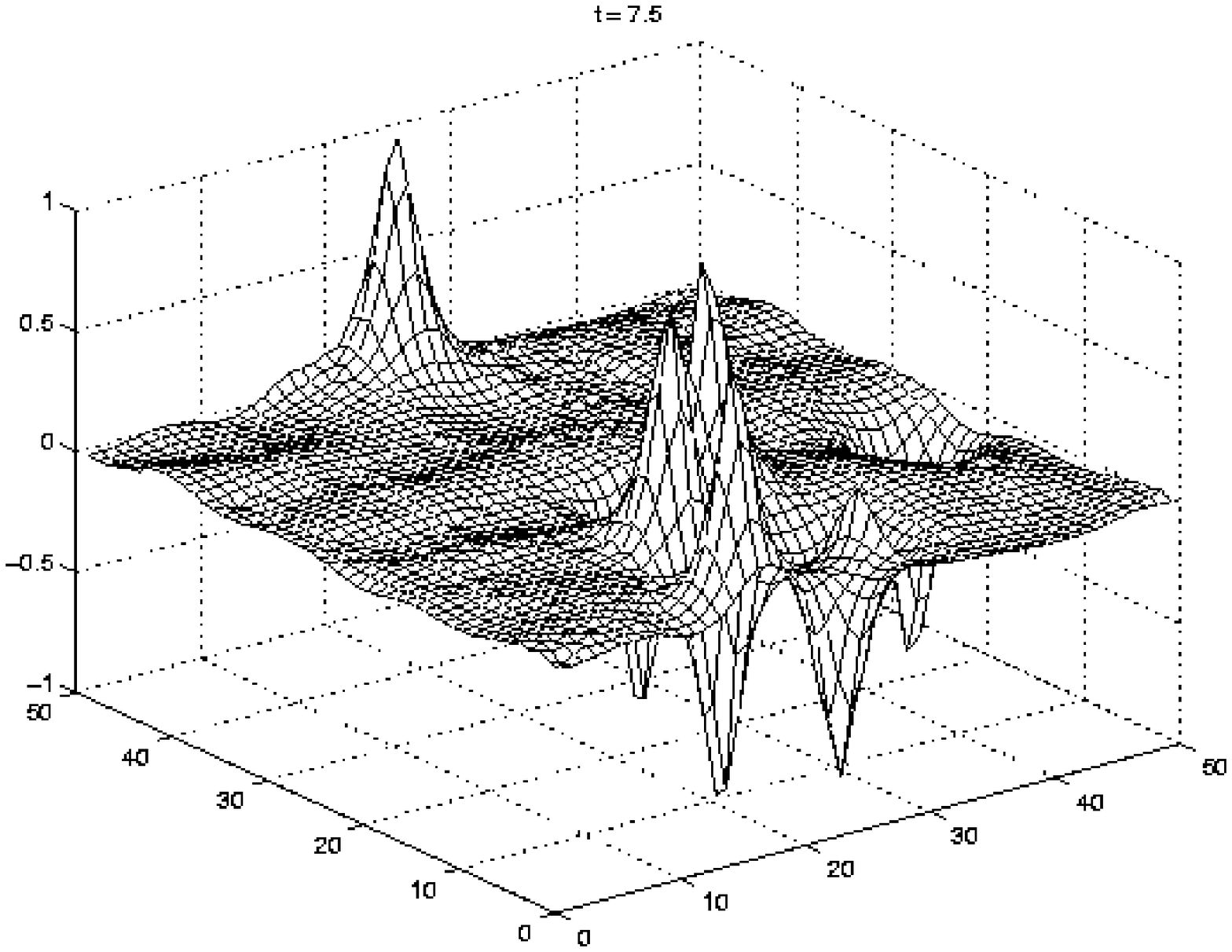,width=7cm} \\
      $ t=6$   &   $ t=7.5$ \\
      \epsfig{file=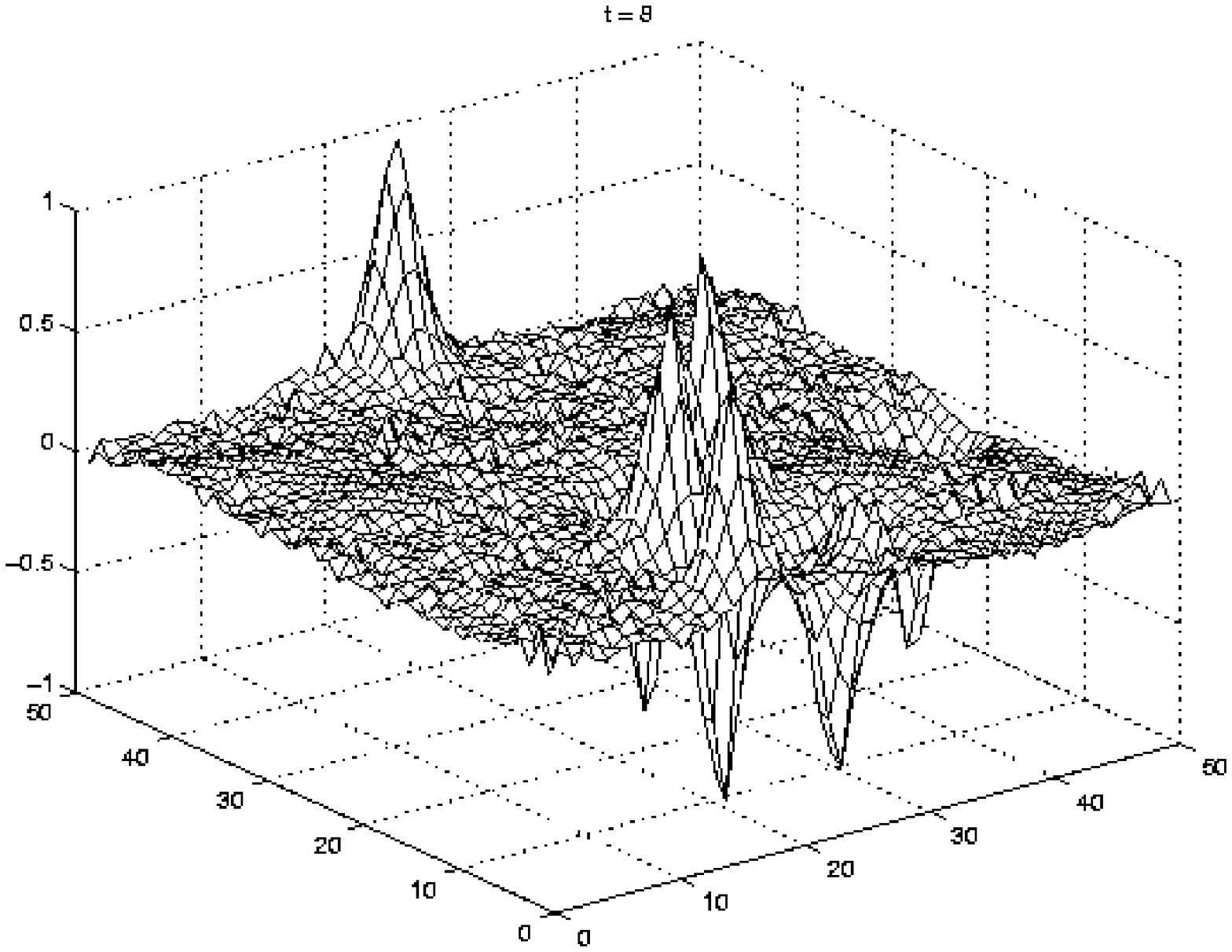,width=7cm}
      & \epsfig{file=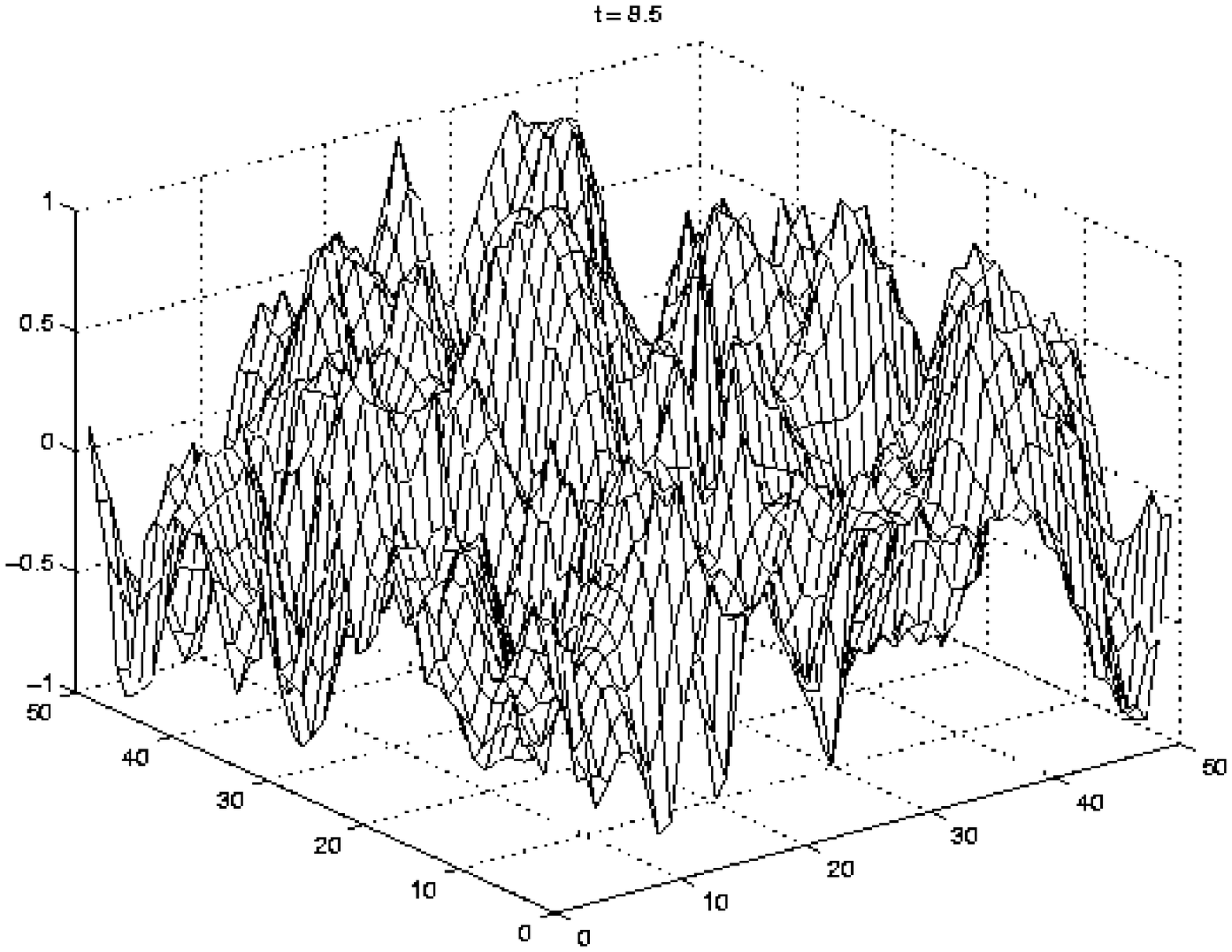,width=7cm} \\
      $ t=8$   &   $ t=8.5$
    \end{tabular}
    \caption{Snapshots of ``wandering vortices'', Example 3.}
    \label{fig:8}
  \end{center}
\end{figure}

\begin{figure}[htbp]
  \begin{center}
    \begin{tabular}{lr}
      \epsfig{file=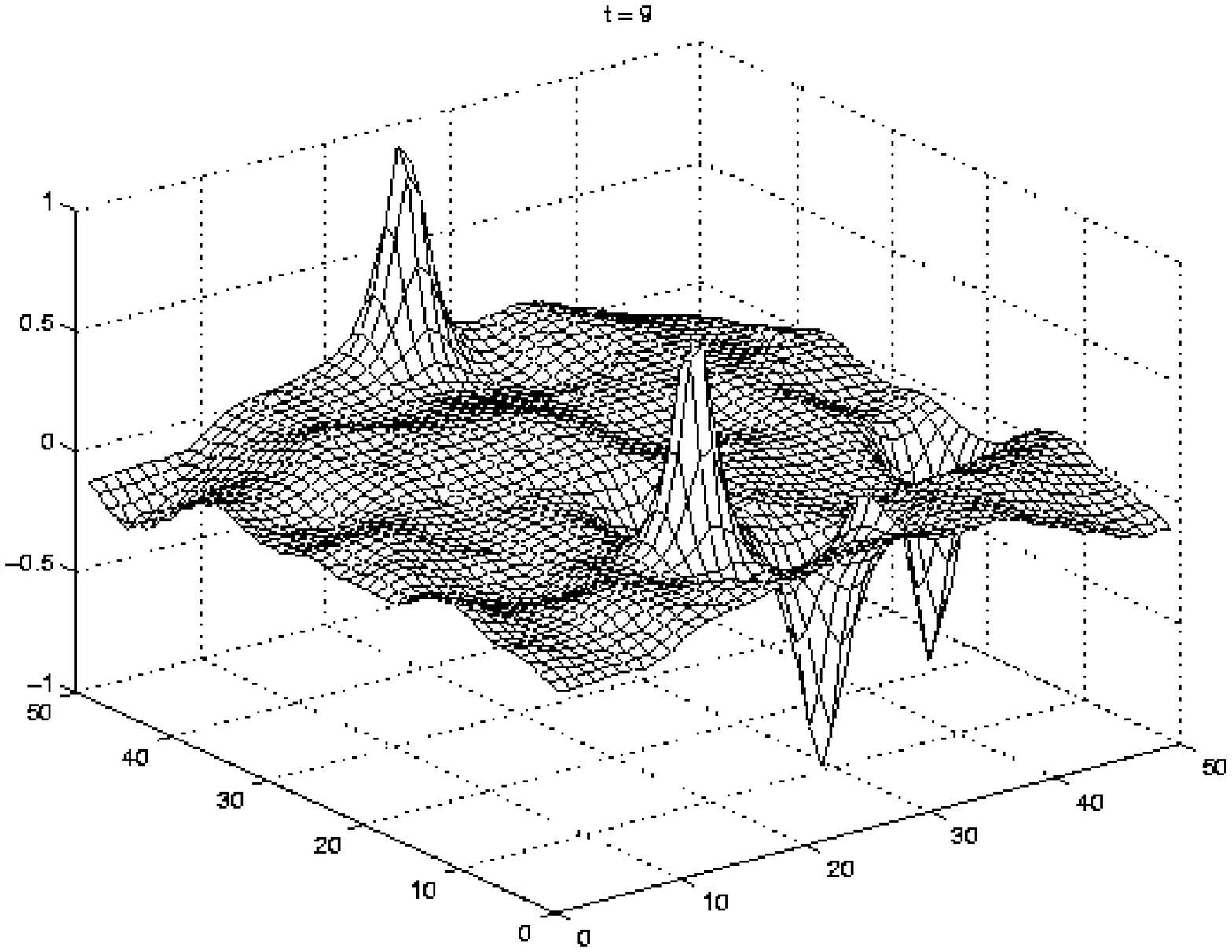,width=7cm}
      & \epsfig{file=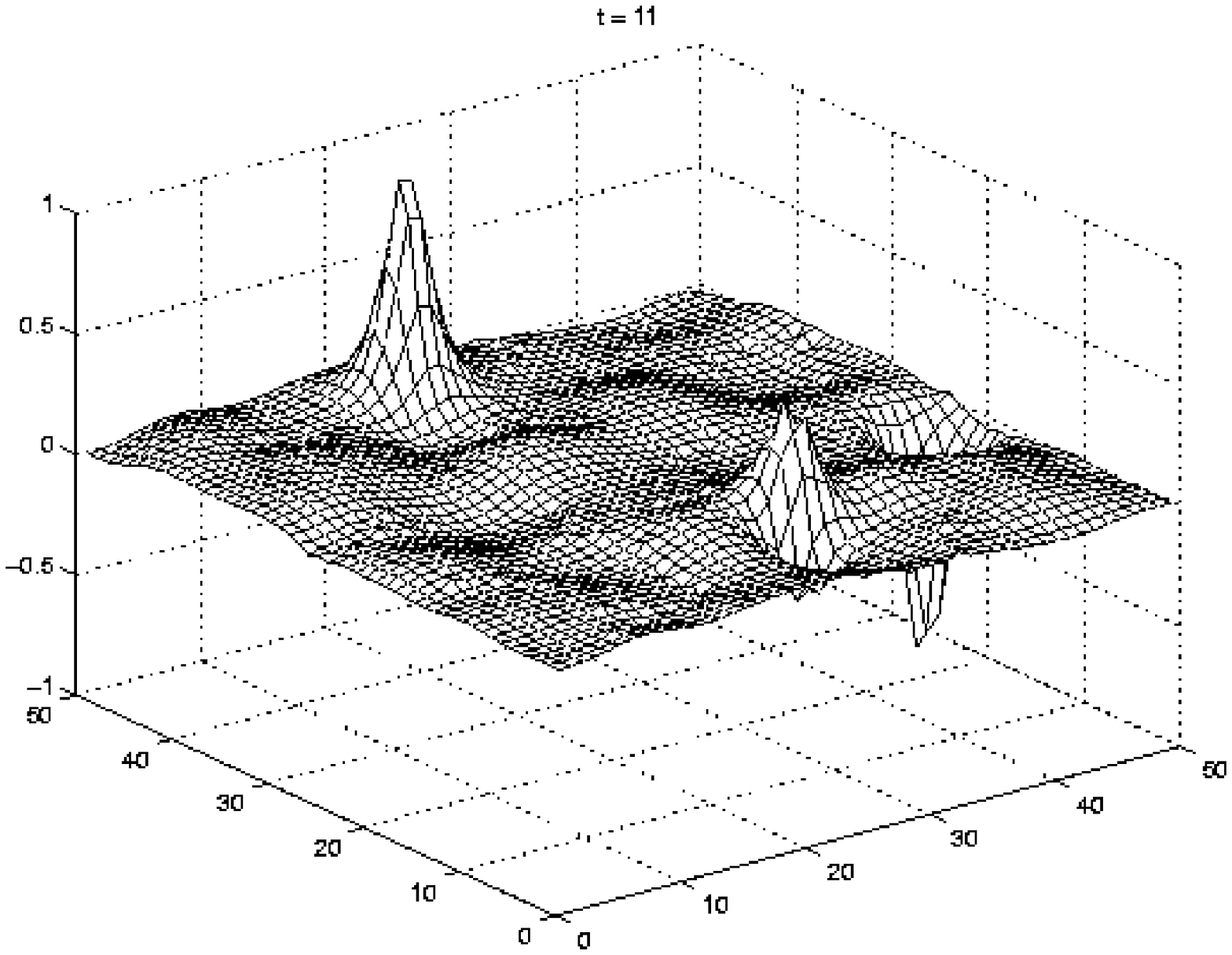,width=7cm} \\
      $ t=9$   &   $ t=11$ \\
      \epsfig{file=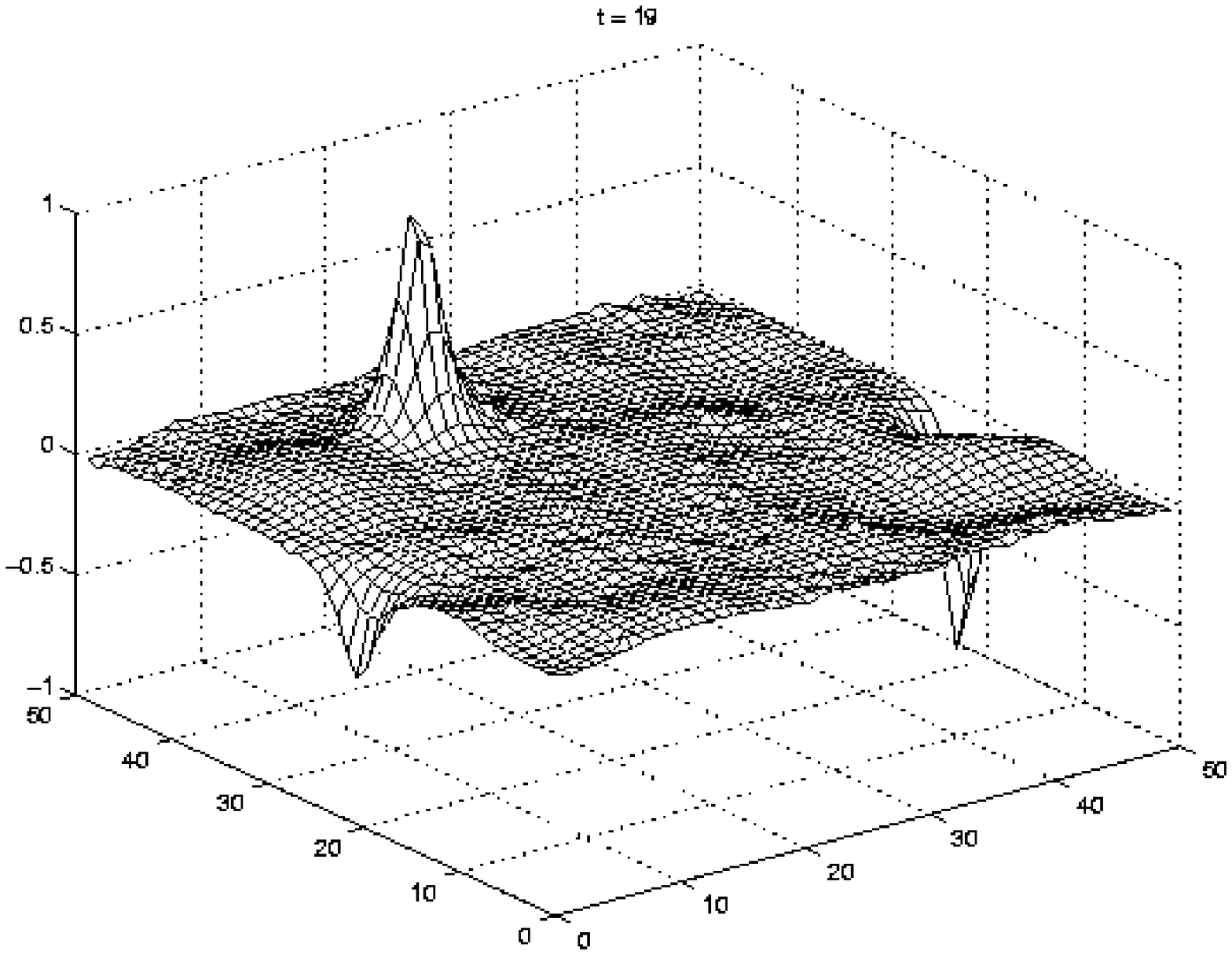,width=7cm}
      & \epsfig{file=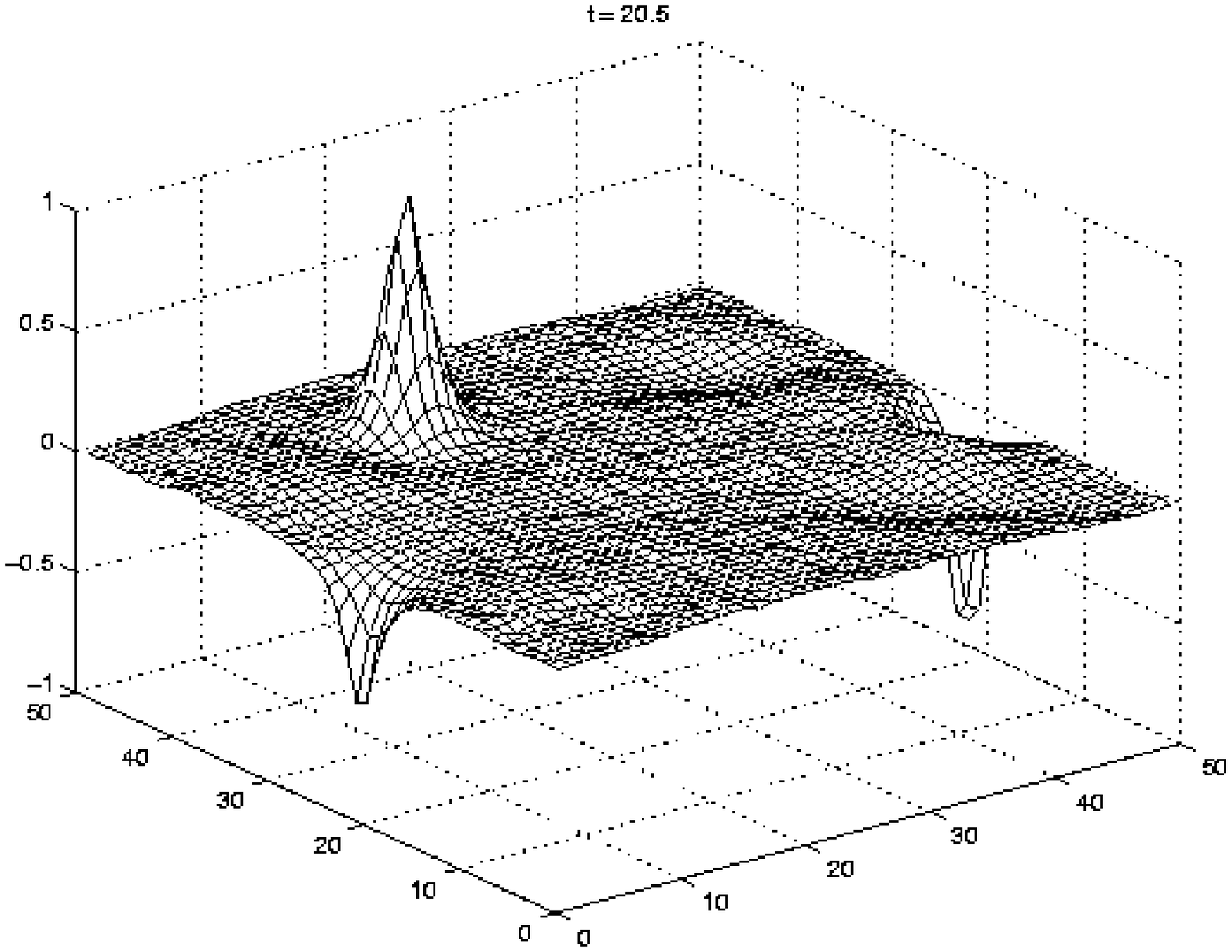,width=7cm} \\
      $ t=19$   &   $ t=20.5$ \\
      \epsfig{file=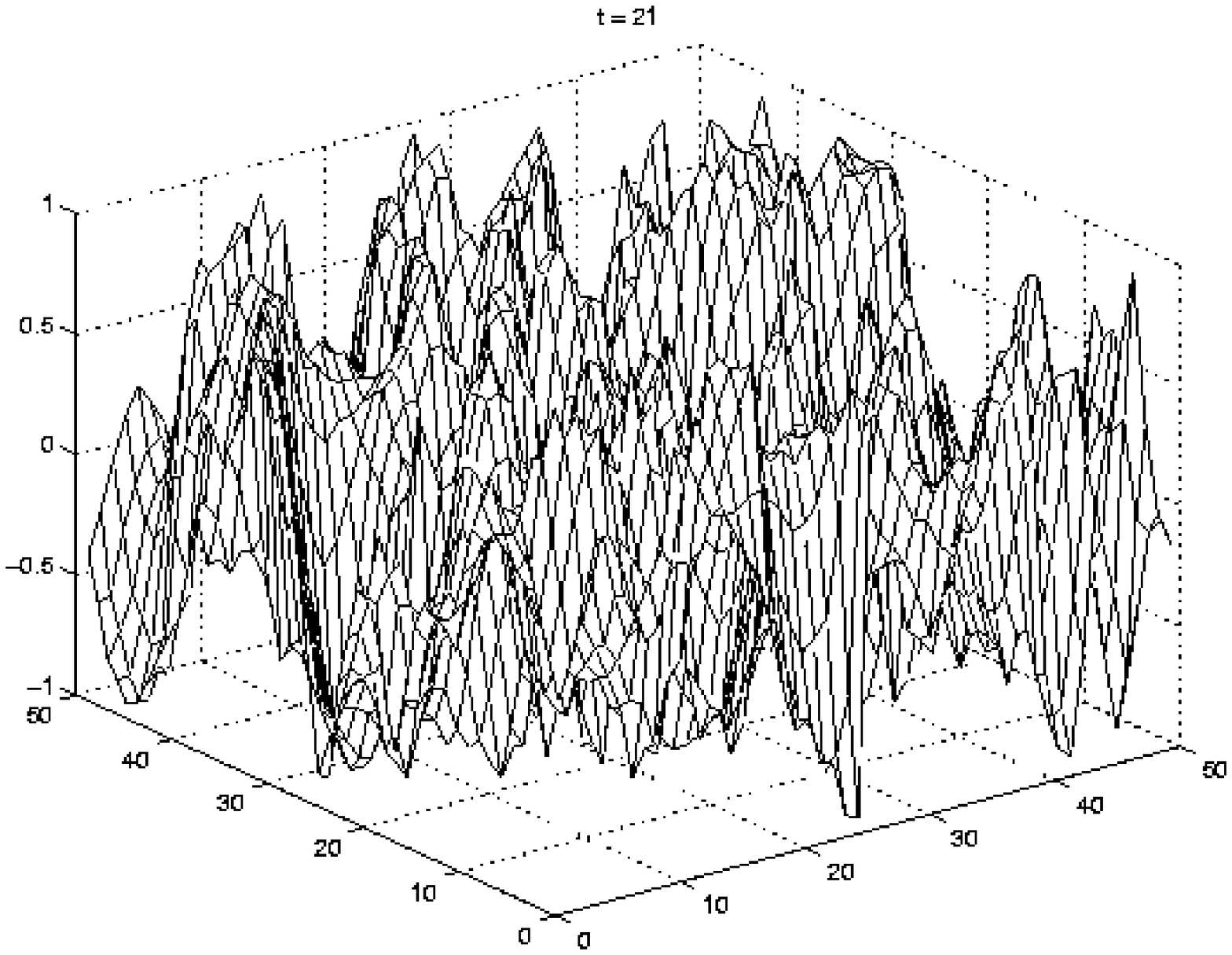,width=7cm}
      & \epsfig{file=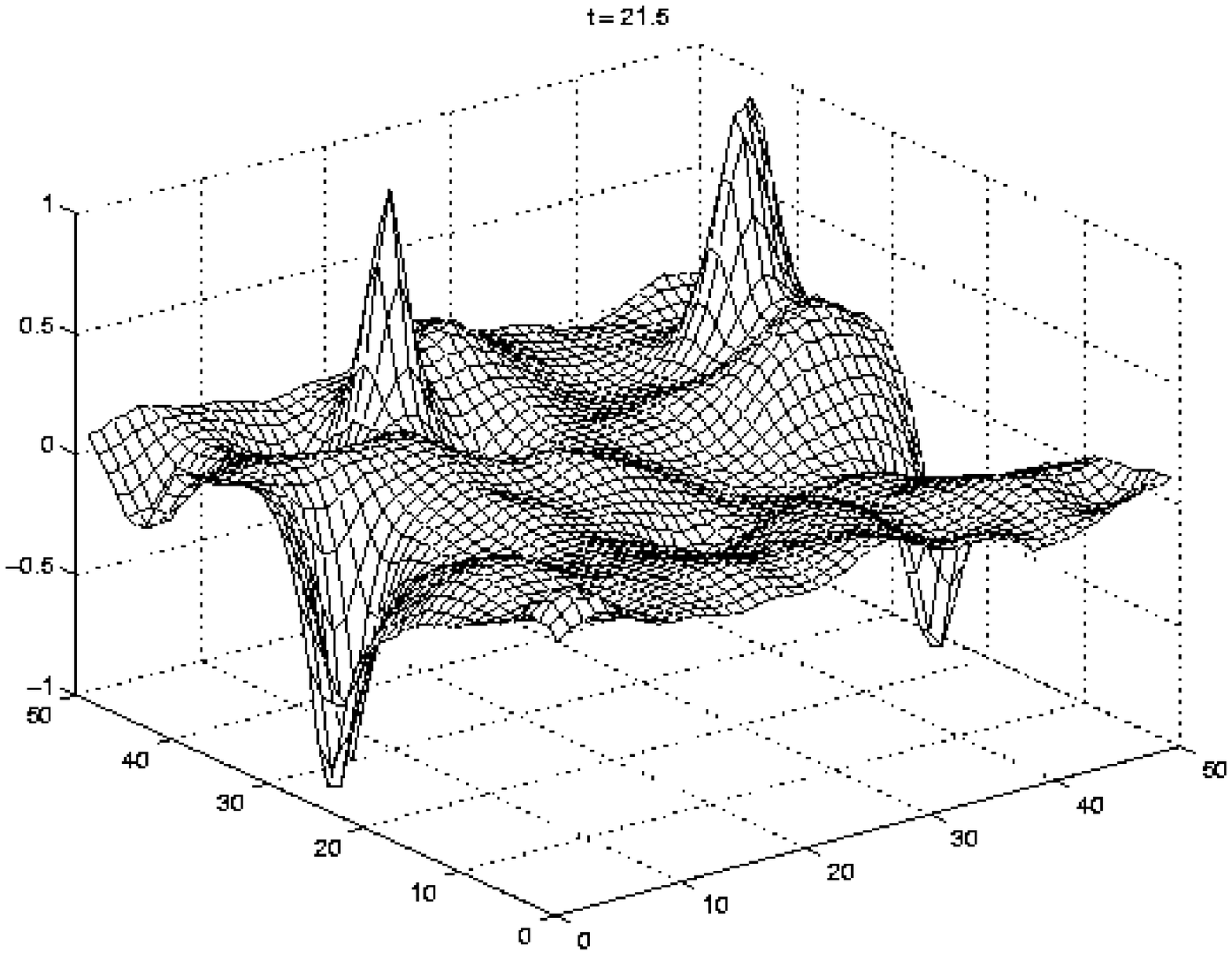,width=7cm} \\
      $ t=21$   &   $ t=21.5$
    \end{tabular}
    \caption{Snapshots of ``wandering vortices'', Example 3.}
    \label{fig:9}
  \end{center}
\end{figure}

{\em Example 3.}  Another example, in Figures \ref{fig:5}-\ref{fig:9}
is a thermostatted system we call ``wandering vortices''.  This
beautiful system has random initial conditions, periodic boundaries,
and parameters $\lambda=0.9$, $T=0.05$ and timestep $\Delta t=0.05$.
Figures \ref{fig:5} and \ref{fig:6} represent the evolution of
$\alpha$, energy, and averages of the energy over different time
windows.  Interestingly, the behavior of $\alpha$ is much wilder than
in the Example 2.  The snapshots in Figures
\ref{fig:8},\ref{fig:9} show the $z-$components of the lattice.  From
random state, the system very quickly forms vortices on smooth
surrounding, which move around for a short time, then look random
again, then vortices again.  Sometimes the vortices died out
completely leaving us just with smooth surface, then reappeared again.

\begin{figure}[htbp]
  \begin{center}
    \begin{tabular}{lr}
      \epsfig{file=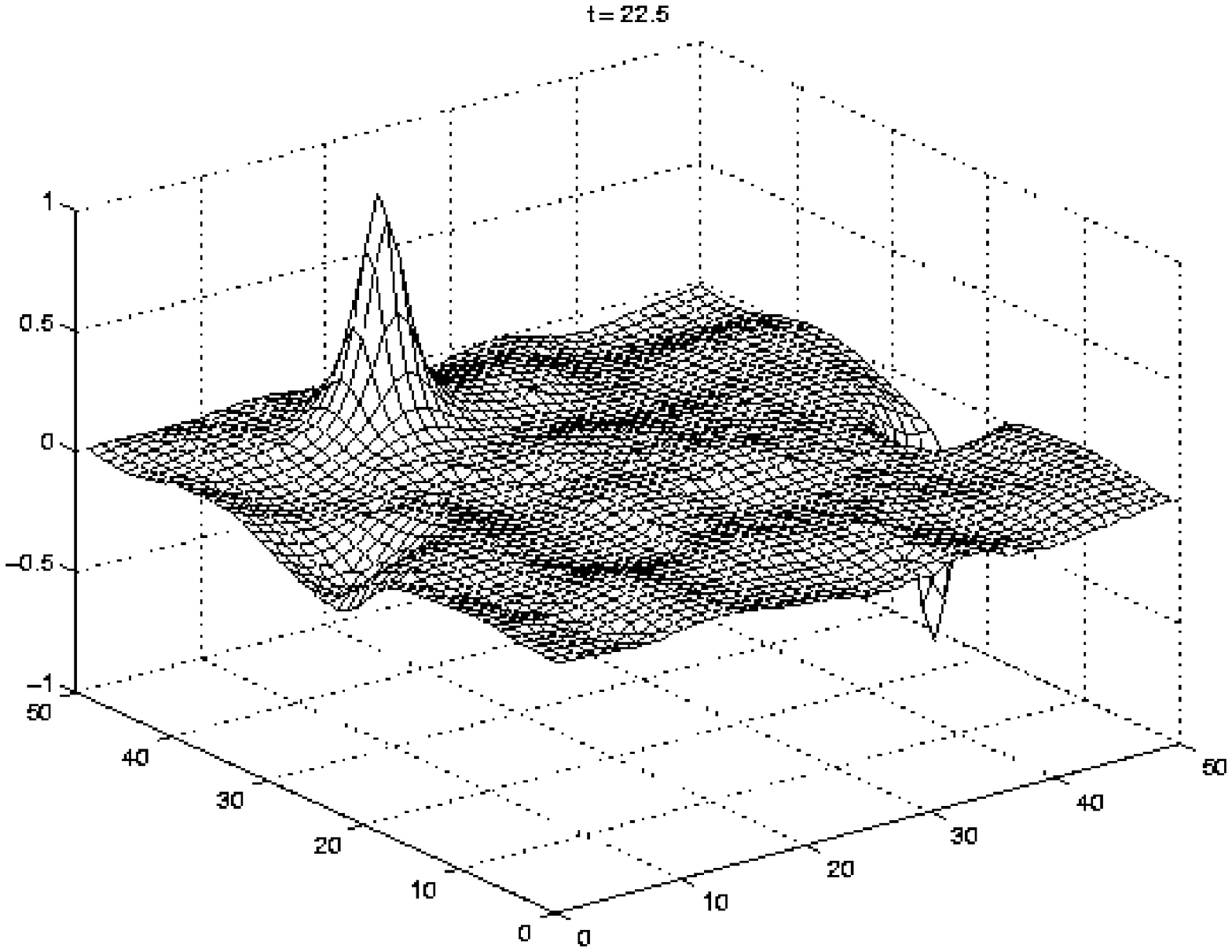,width=7cm}
      & \epsfig{file=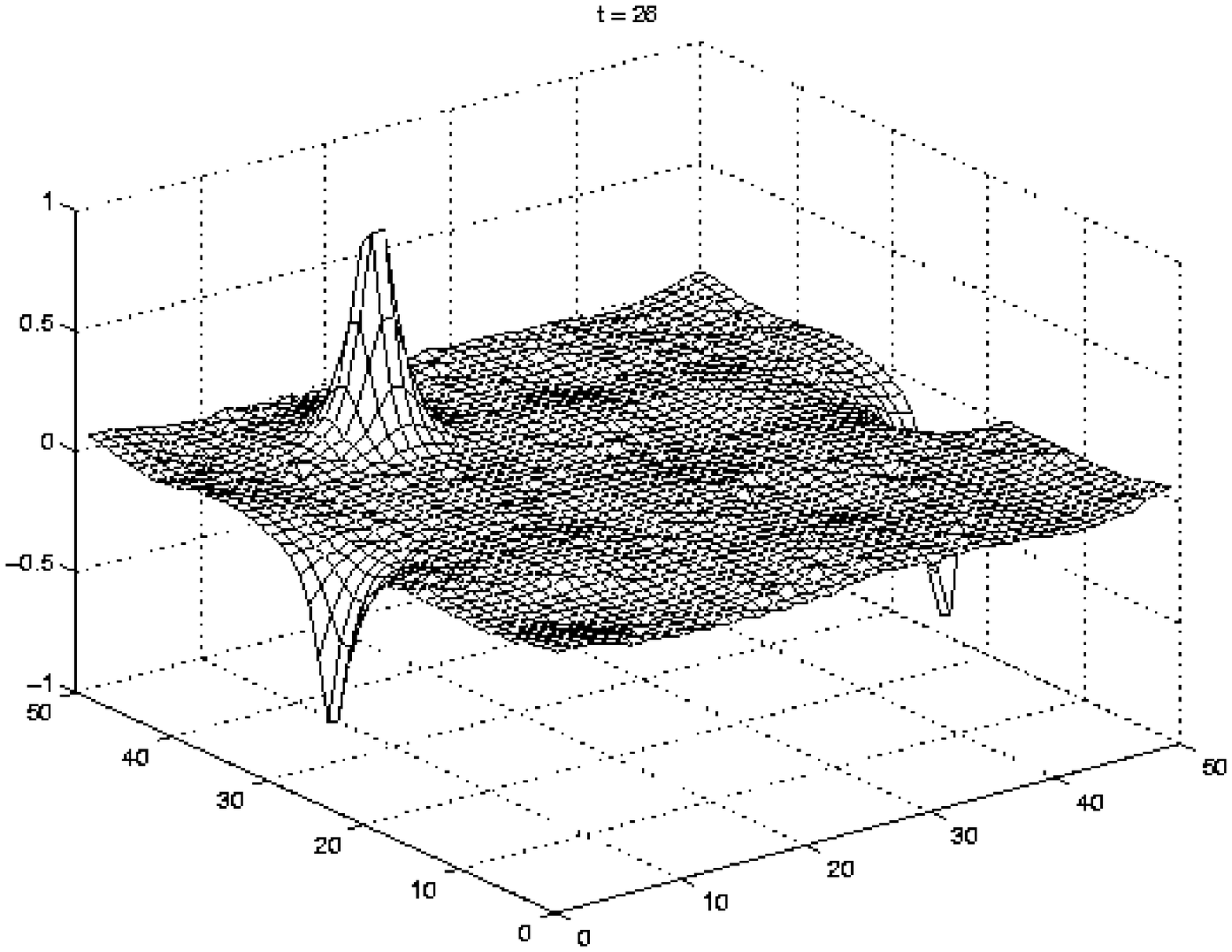,width=7cm} \\
      $ t=22.5$   &   $ t=26$ \\
      \epsfig{file=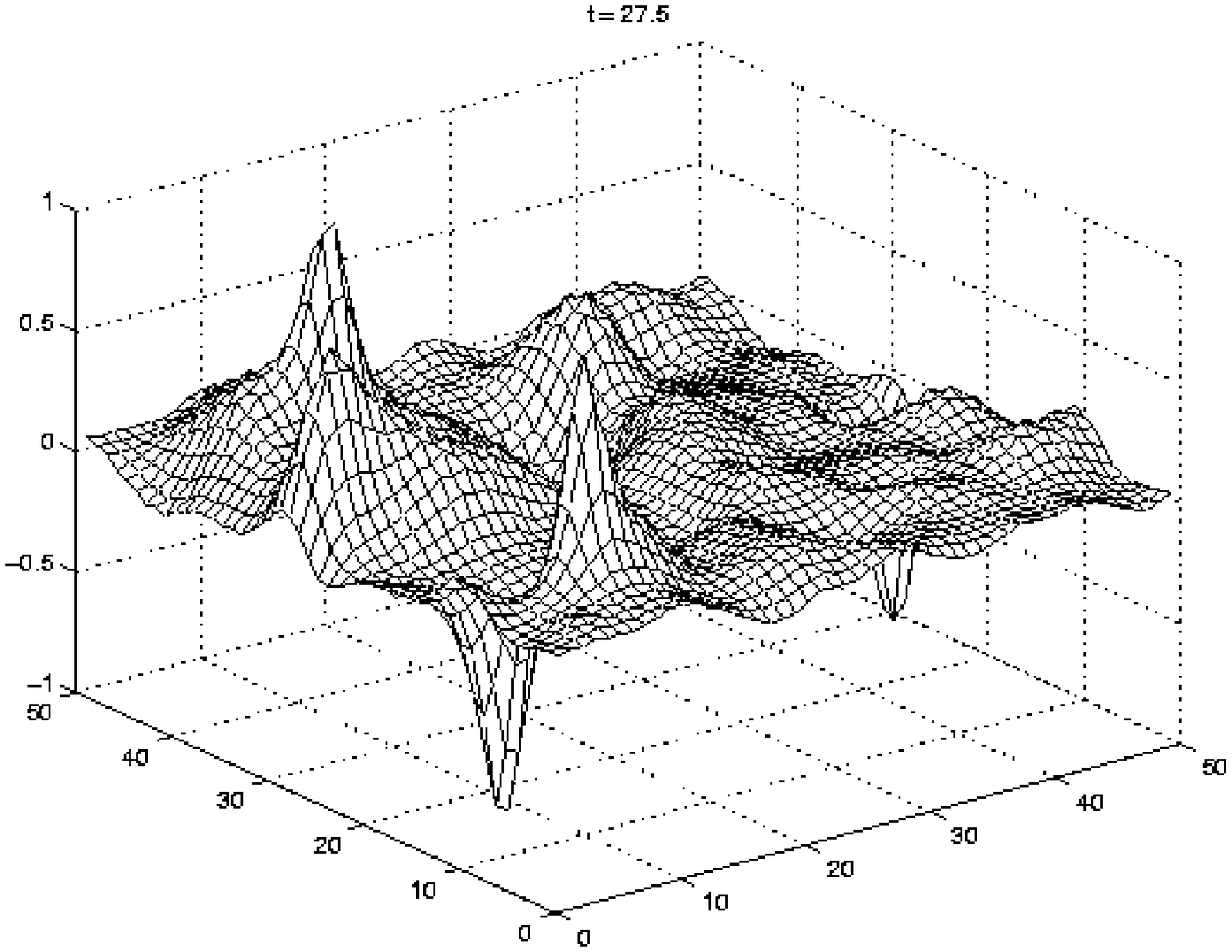,width=7cm}
      & \epsfig{file=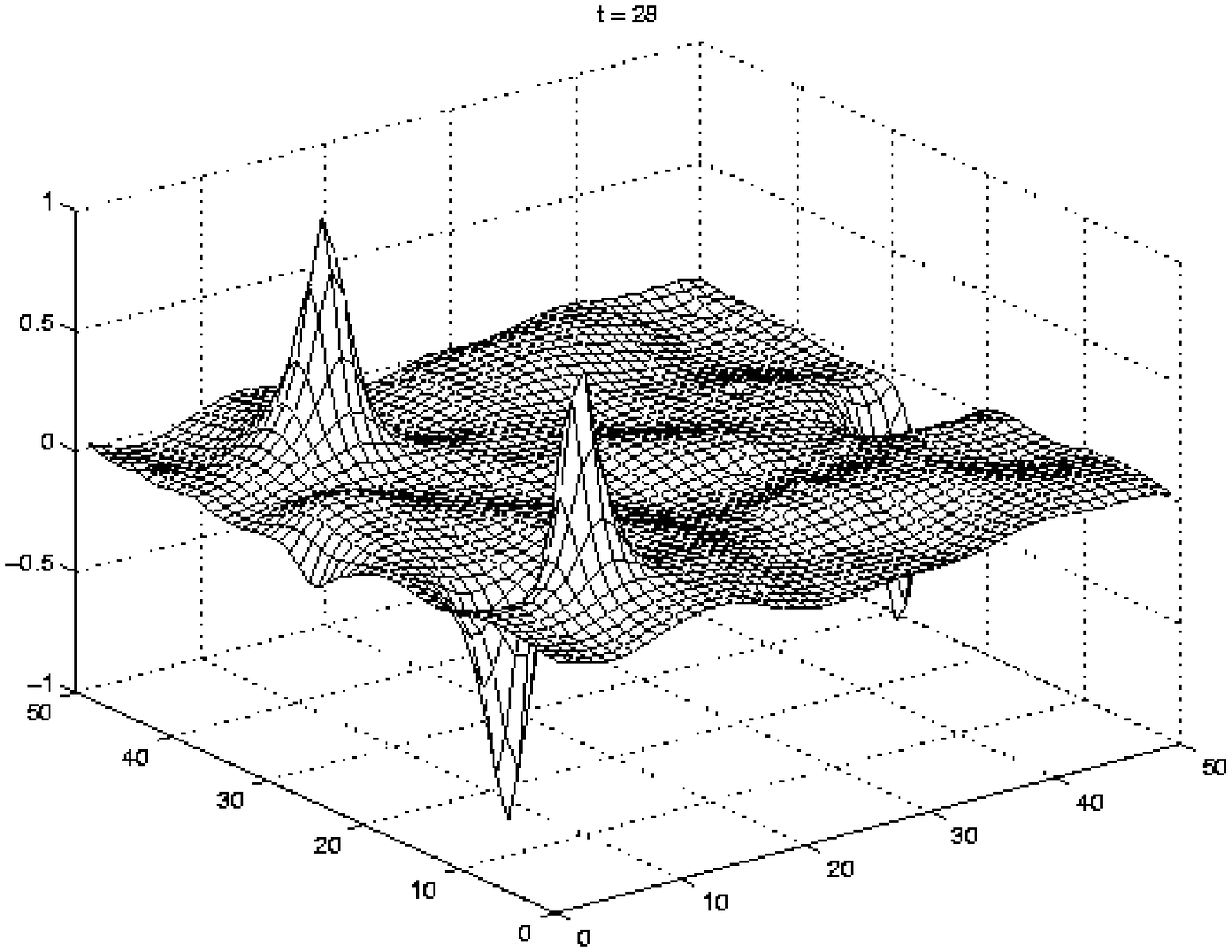,width=7cm} \\
      $ t=27.5$   &   $ t=28$ \\
      \epsfig{file=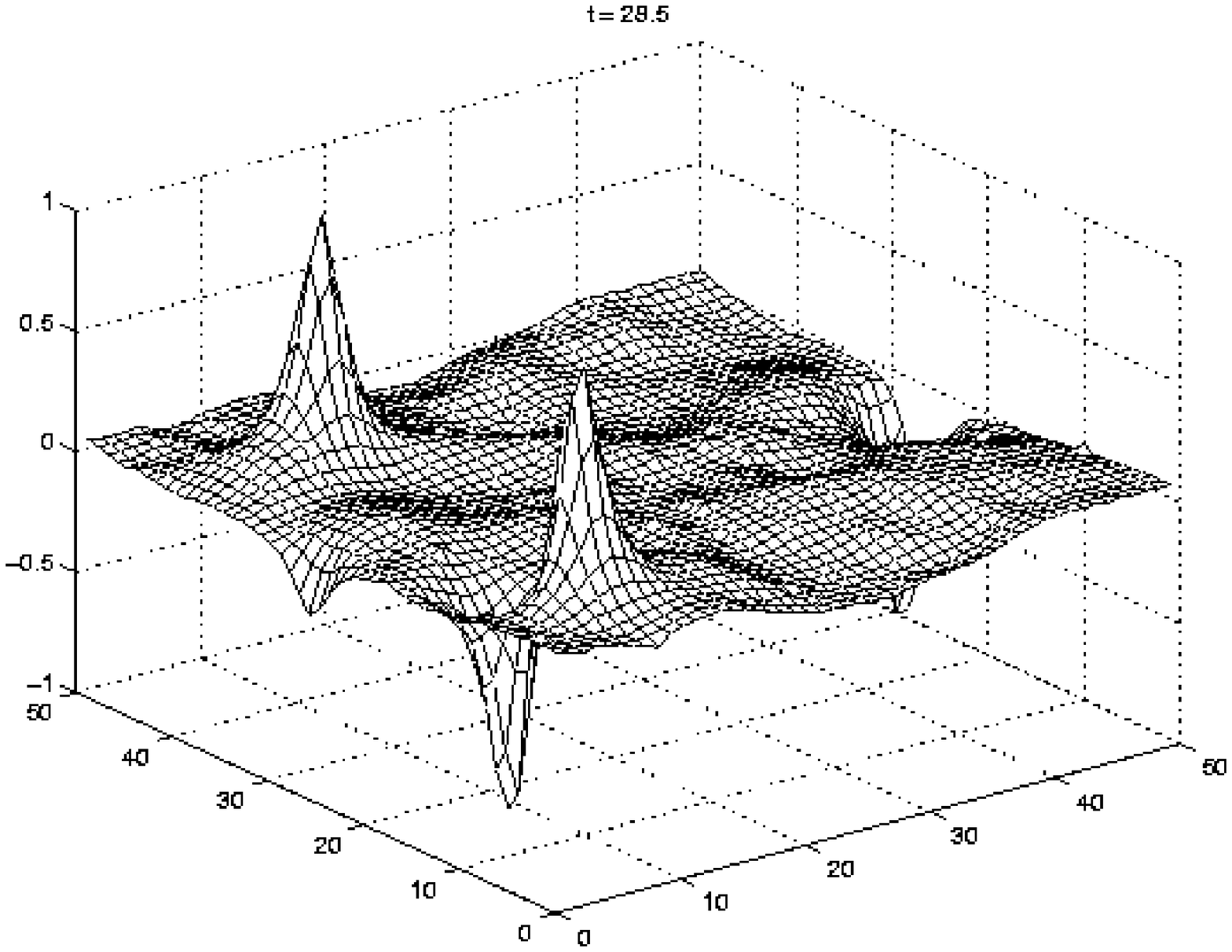,width=7cm}
      & \epsfig{file=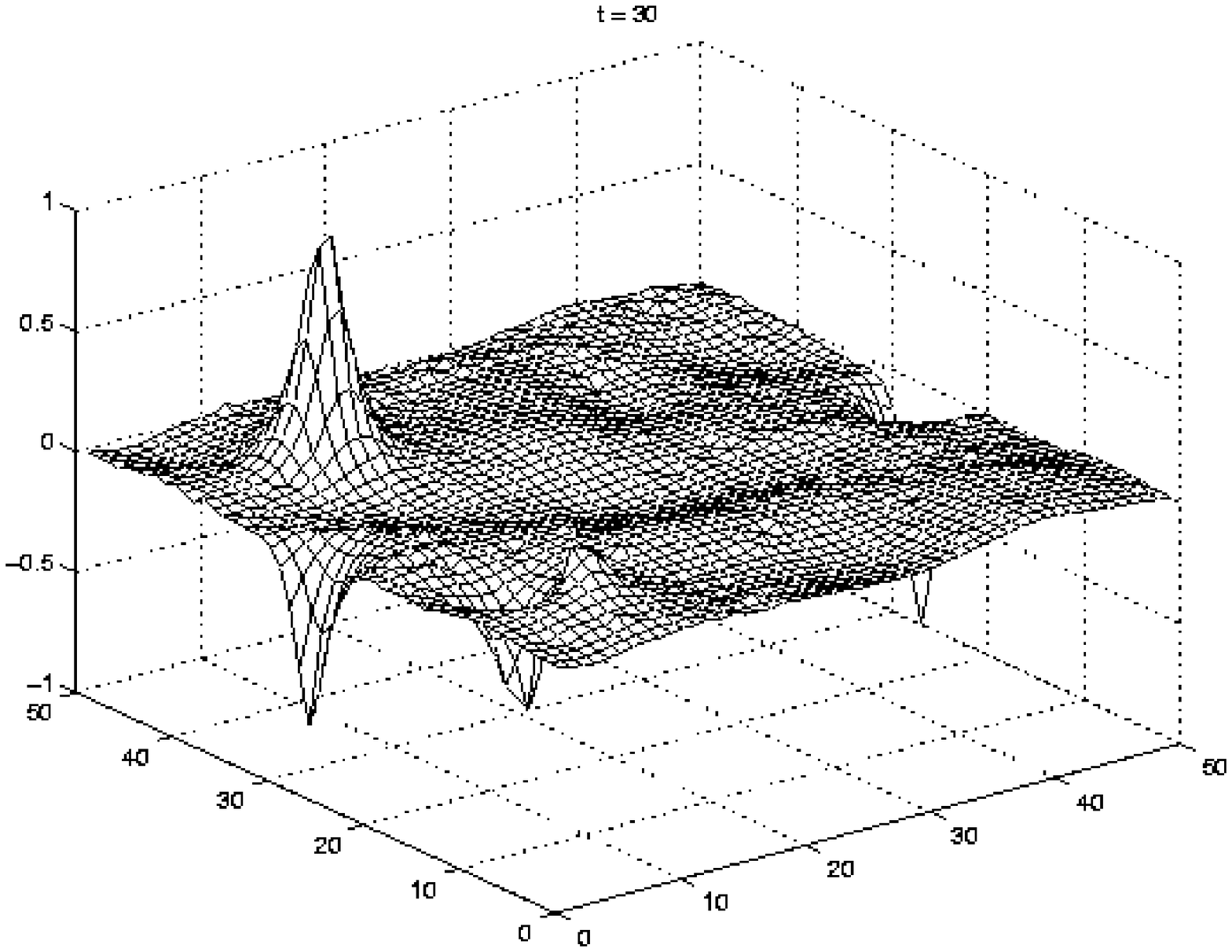,width=7cm} \\
      $ t=28.5$   &   $ t=30$
    \end{tabular}
    \caption{Snapshots of ``wandering vortices'', Example 3.}
    \label{fig:10}
  \end{center}
\end{figure}

\subsubsection*{Comparison to RK4 with projection}

For comparison we implemented the classical Runge-Kutta $4^{th}$ order
method (RK4) with projection: after every step we normalize
\[
  z_{ij,new} := z_{ij,RK}/|z_{ij,RK}|,
\]
where $z_{ij,RK}$ denotes the result of RK4 step.  At very small
timesteps for which the RK4 method could successfully integrate the
problem, it was slightly faster than splitting, but the RK4 method
became rapidly unstable as the stepsize and/or anisotropy were
increased.  The splitting method was able to handle large anisotropies
($\lambda=3$) and step sizes ($\Delta t=0.3$).  However, we did not
seek the limits of our splitting method.  The values $\lambda=3$ and
$\Delta t=0.3$ indicate the superior stability well enough at this
stage.

We tested this projected RK4 on Examples 2 and 3.  We kept the other
parameters intact but changed the step size.  As a sign of failure, we
stopped the computation when the code started to produce infinities.
In Example 2 the maximal timestep was 0.01, while in Example 3 $\Delta
t_{max} = 0.015$.  The results are summarized in table \ref{tab:1}.
\begin{table}[htbp]
  \begin{center}
    \begin{tabular}[t]{|c|c|c|}
      \hline
      Example 2 & $\Delta t$ & RK4 survived until $t=$ \\
      & 0.020 &  0.060 \\
      & 0.015 &  0.075 \\
      & 0.012 &  0.084 \\
      & 0.011 &  0.099 \\
      & 0.0105 &  0.116 \\
      & 0.0102 &  0.235 \\
      & 0.010  &  $>200$ \\
      \hline
      Example 3 & $\Delta t$ & RK4 survived until $t=$ \\
      & 0.020 &  0.080 \\
      & 0.019 &  0.095 \\
      & 0.017 &  0.102 \\
      & 0.016 &  0.112 \\
      & 0.015 &  $>150$ \\
      \hline
    \end{tabular}
    \caption{Failure of projected RK4}
    \label{tab:1}
  \end{center}
\end{table}

\section{Discussion}
\label{sec:discuss}

In this paper we have developed and tested a geometric integrator for
a semi-discretized Landau-Lifshitz-Gilbert (LLG) equation which
includes a nonlinear dissipative term, as well as for a more
complicated thermostatted model, following the approach of Bulgac and
Kusnetsov.

The integrator for the dissipated system is shown to have a dissipative
property.  However, it is difficult to compare since we do not know
the exact continuous solution.  LLG is currently a very active topic
of research, see more details in the introduction.  However, it seems
that so far there has not been developed a geometric integrator for a
thermostatted system.

Trying to simulate the thermostatted systems with projected RK4
revealed both the features of a stiff ODE and features of a
conservative system.  The combination is extremely difficult for
standard form numerical methods.  The key feature of our splitting
method is that it is constructed from composition of building blocks
that simulate each of the two components of the system correctly.

Simulation with our new thermostatted method has revealed interesting
phenomena: slowly creeping boundaries, or wandering vortices, both of
which appear from random initial conditions.  
Our informal term  ``wandering vortices'' in Example 3 refers not to
certain particular vortices that survive throughout the whole simulation, 
but to a situation where
we have two or more vortices which wander for a while, then violently crash
and form new vortices.
Intermediate states
include ``quasi-chaotic'' state, an informal term by which we mean a
state that suddenly appears and looks random but is not, since it
renders immediately back to (almost) the same smooth motion.

The RK4 method is less stable.   The stepsize restriction is an order
of magnitude smaller compared to our splitting method.    This is
evidence of stiffness in the ODEs, and a better choice might seem to be
a stiff solver on this account, but if one uses a stiff solver the
result would be poor resolution of the conservative evolution which is
also an important component of the dynamics of the system.   The best
compromise is therefore a composition scheme, such as that outlined
here, which separately  and appropriately resolves each term of the
system.

We anticipate that this work will stimulate further research in the
development of thermostatted numerical methods for systems with
complicated nonlinear structure."

\noindent{\bf Acknowledgement.} 
The authors are grateful to Jason Frank both for valuable comments and
providing the code of \cite{fr-hu-le97:GIC}.  The author TA was
supported by the Academy of Finland.  BL was supported by the
Engineering and Physical Sciences Research Council, grant
GR/R03259/01.

\section*{Appendix: Lie-Poisson Canonical Sampling Technique}
Consider first the Poisson rigid body system consisting of a Hamiltonian $H(z)$ together with structure matrix $J(z)$ admitting
Casimir $|z|$.   Based partly on \cite{la-le03:GDT} we construct an augmented Hamiltonian with additional variables
$\sigma$, $\pi_{\sigma}$, $\theta$, $\pi_\theta$:
\[
\tilde{H} = H(Q(\theta)z) + \pi_\sigma^2/2\mu + kT \ln \sigma + G(\theta, \pi_\theta, \pi_\sigma).
\]
where $Q$ is an orthogonal matrix depending on parameter(s) $\theta$.    The Lie-Poisson structure is just the
rigid body Poisson structure augmented by the canonical structure for the augmenting variables. Under assumption of 
ergodicity, and some very mild technical conditions similar to those obtained in \cite{la-le03:GDT}, this Hamiltonian can be shown to provide canonical sampling from microcanonical trajectories i.e. 
\[
\int \int \int \int f(Q(\theta) z) \delta [ \tilde{H} - E ] d\sigma d \theta d \pi_\sigma d\pi_\theta = 
f(z) \exp ( -\frac{1}{kT} H(z) ),
\]
with preservation of the Casimir due to $Q$ being orthogonal.
In order to obtain ergodicity, the "bath Hamiltonian" $G$ should be sufficiently complicated.

In the case of a spin system, $H(z_1,z_2,\ldots,z_N)$,  we may introduce a separate unit 3-vector $\theta_i$ for each spin vector.  Then the
Hamiltonian
\[
\tilde{H} = H(Q(\theta_1)z_1,Q(\theta_2)z_2,\ldots,Q(\theta_N)z_N) + \pi_{\sigma}^2/2\mu + kT \ln \sigma + G(\theta, \pi_{\theta}, \pi_{\sigma}),
\]
will enable canonical sampling.   For example $Q(u)$ can be a Householder transformation, 
\[
Q(u) = I - 2 \frac{u u^T}{u^T u},
\]
and the bath Hamiltonian $G$ can describe a coupled system of spherical pendula involving in some nontrivial way the parameter $\pi_{\sigma}$.

\end{document}